\documentclass[preprint,authoryear]{elsarticle}

\usepackage{graphicx}

\usepackage{amssymb, amsmath}
\usepackage{subfigure}

\usepackage[breaklinks=true,colorlinks=true]{hyperref}

\usepackage{booktabs}
\usepackage{multirow}

\newcommand{\spr}[1]{\ensuremath{^\textrm{#1}}}

\renewcommand{\eqref}[1]{Eq.~(\ref{#1})}

\newcommand{\ssun}{\textsc{sun}}
\newcommand{\A}{\textsc{a}}
\newcommand{\B}{\textsc{b}}
\newcommand{\C}{\textsc{c}}

\newcommand{\h}{\textsc{h}}

\newcommand{\La}{\textsc{l}}
\newcommand{\M}{\textsc{m}}

\newcommand{\R}{\textsc{r}}
\newcommand{\s}{\textsc{s}}
\newcommand{\T}{\textsc{t}}

\newcommand{\AU}{\textsc{au}}
\newcommand{\LIM}{\textsc{lim}}

\journal{Advances in Space Research}

\graphicspath{{gfx/}}

\begin{document}

\begin{frontmatter}



\title{Design of a Formation of Solar Pumped Lasers for Asteroid Deflection}


\author{Massimiliano Vasile\corref{cor}\fnref{footnote2}}
\address{Dept.\ of Mechanical \& Aerospace Engineering, University of Strathclyde, 75 Montrose Street, G1~1XJ, Glasgow, UK, Ph. +44 (0)141 548 2083, Fax +44 (0)141 552 5105}
\cortext[cor]{Corresponding author}
\fntext[footnote2]{Reader in Space System Engineering}
\ead{massimiliano.vasile@strath.ac.uk}


\author{Christie Alisa Maddock\fnref{footnote3}}
\address{Dept.\ of Mechanical \& Aerospace Engineering, University of Strathclyde, 75 Montrose Street, G1~1XJ, Glasgow, UK}
\fntext[footnote3]{Research Fellow, Centre for Future Air-Space Transport Technology}
\ead{christie.maddock@strath.ac.uk}

\begin{abstract}

This paper presents the design of a multi-spacecraft system for the deflection of asteroids. Each spacecraft is equipped with a fibre laser and a solar concentrator. The laser induces the sublimation of a portion of the surface of the asteroid, and the resultant jet of gas and debris thrusts the asteroid off its natural course. The main idea is to have a formation of spacecraft flying in the proximity of the asteroid with all the spacecraft beaming to the same location to achieve the required deflection thrust. The paper presents the design of the formation orbits and the multi-objective optimisation of the formation in order to minimise the total mass in space and maximise the deflection of the asteroid. The paper demonstrates how significant deflections can be obtained with relatively small sized, easy-to-control spacecraft.

\end{abstract}

\begin{keyword}
Asteroid deflection \sep Laser ablation \sep NEO
\end{keyword}

\end{frontmatter}

\parindent=0.5 cm

\section{Introduction}\label{intro}

\renewcommand*{\thefootnote}{\fnsymbol{footnote}}  

Near Earth Objects (NEO), the majority of which are asteroids, are defined as minor celestial objects with a perihelion less than 1.3~AU and an aphelion greater than 0.983~AU. A subclass of these, deemed Potentially Hazardous Asteroids (PHA), are defined as those with a Minimum Orbital Intersection Distance (MOID) from the Earth's orbit less than or equal to 0.05 AU and a diameter larger than 150~m (equivalent to an absolute magnitude of 22.0 or less). As of March 2012, 8758 NEO's have been detected \citep{mpc}; of those, 840 are estimated to have an effective diameter larger than 1~km\footnote{An asteroid with an effective diameter equal to or greater than 1 km is defined here to be any NEA with an absolute brightness or magnitude $H \le 17.75$, as per \citet{Stuart2003}.}, and 1298 are categorised as potentially hazardous. Impacts from asteroids over 1 km in diameter are expected to release over $10^5$ megatons of energy with global consequences for our planet \citep{stokes2003}, while those with an average diameter of 100~m can are expected to release over $10^2$ megatons of energy potentially causing significant tsunamis and/or land destruction of a large city \citep{Toon1997}. It is estimated that there are between 30,000--300,000 NEO's with diameters around 100~m, meaning a large number of NEO's are still undetected.

A quantitative comparison of the various options for NEO deflection was conducted by \citet{Colombo2006,Sanchez2007}. Examining the results of the comparison, one of the more interesting methods employed solar sublimation to actively deviate the orbit of the asteroid. The original concept, initially evaluated by \citet{Melosh1994}, and later assessed by \citet{Kahle2006}, envisioned a single large reflector; this idea was expanded to a formation of spacecraft orbiting in the vicinity of the NEO, each equipped with a smaller concentrator assembly capable of focusing the solar power at a distance around 1 km and greater \citep{Maddock2007}. This concept addressed the proper placement of the concentrators in close proximity to the asteroid while avoiding the plume impingement and provided a redundant and scalable solution. However, the contamination of the optics still posed a significant limitation as demonstrated by \citet{vasile2010}. In the same paper, the authors demonstrated that the combined effect of solar pressure and enhanced Yarkovsky effect could lead to miss (or deflection) distances of a few hundred to a thousand kilometres over eight years of deflection time. However, this deflection is orders of magnitude lower that the one achievable with a prolonged sublimation of the surface.

A possible solution is to use a collimating device that would allow for larger operational distances and protection of the optics. This paper presents an asteroid deflection method based on a formation of spacecraft each equipped with solar pumped lasers. The use of lasers has already proposed by several authors, although always in conjunction with a nuclear power source \citep{Phipps1992,Phipps1997,Park2005}. Extensive studies on the dynamics of the deflection with high power lasers were proposed by \citet{Park2005} envisaging a single spacecraft with a MW laser. This paper proposes a different solution with a formation of smaller spacecraft, each supporting a kW laser system indirectly pumped by the Sun.

The paper starts with a simple sublimation model that is used to compute the deflection force. The orbits of the spacecraft formation are then designed by solving a multi-objective optimisation problem that yields an optimal compromise between distance from the target and impingement with the plume of debris. A Lyapunov controller is proposed to maintain the spacecraft in formation along the desired proximal orbit. A second multi-objective optimisation problem is then solved to compute a different type of controlled formation orbits in which the shape of the orbit is predefined. Finally, the number and size of the spacecraft is optimised to yield the maximum possible deflection.  

\section{Deflection Model}\label{sec:1}

The orbital properties of Near Earth Asteroids (NEA) can be grouped into four general categories based on the semi-major axis $a$ of the orbit, radius of apoapsis $r_a$ and/or radius of periapsis $r_p$, described as follows \citep{nasa_neo}:
\begin{description}
  \item[Atens] Earth-crossing asteroids with semi-major axes smaller than Earth (named after asteroid 2062 Aten), where $a < 1$ AU, $r_a \geq 0.983$ AU.
  \item[Apollos] Earth-crossing asteroids with semi-major axes larger than Earth (named after asteroid 1862 Apollo), where $a \geq 1$ AU, $r_p \leq 1.0167$ AU.
  \item[Amors] Earth-approaching asteroids with orbits exterior to Earth's but interior to Mars (named after asteroid 1221 Amor), where $a > 1$ AU, $1.0167$ AU $< r_p \leq 1.3$ AU.
  \item[Atiras] Near Earth Asteroids whose orbits are contained entirely with the orbit of the Earth (named after asteroid 163693 Atira), where $a < 1$ AU, $r_p<0.983$ AU.
\end{description}
Apollo is the largest class (approximately 4100 NEA's) followed by Amors (approximately 3400 NEA's). 
The asteroid Apophis 99942, part of the Apollos class, is taken as a test case with a relatively low aphelion such that enough solar power can be harvested.

While circular, or near-circular, orbits offer a more constant level of solar radiation, as suggested by \citet{Vasile2008b,vasile2010} if the mirrors have variable optics, i.e., the focal point can be changed, a constant power density can be achieved for asteroids on elliptical orbits or when the level of solar power available is low.
In terms of altering an orbit, thrusting at the perihelion of elliptical orbits maximises the change in semi-major axis (and therefore the miss distance). Even if the level of solar radiation available is not sufficient to induce sublimation at aphelion, a deflection can still be achieved, as will be demonstrated in this paper.

The other benefit of basing the test case on the Apophis asteroid is its popularity in scientific literature due to the initial, relatively high impact level (2.7\% chance of impacting the Earth in 2029) it was given when it was first observed in 2004. While further tracking data has reduced the threat level, ruling out the possibility of an impact in 2029 but leaving a non-zero impact probability for the 2036 and 2037 encounters, the asteroid Apophis remains a popular reference example. Note, in fact, that although Apophis is not necessary a typical case, the interest here is to examine the effectiveness of a fractionated laser ablation system applied to the deflection of an S-type asteroid, belonging to a given size range, on a moderately eccentric orbit (although also an extension to highly eccentric orbits will be demonstrated) of which Apophis is an example.

Table~\ref{t:apophis} gives the orbital and physical data of the asteroid used in this study. The asteroid shape was assumed to be tri-axial ellipsoidal,
\begin{equation}
a_\ell =  \sqrt{2} d_\A \qquad b_\ell = d_\A \qquad c_\ell = \frac{d_\A}{\sqrt{2}}
\end{equation}
where $a_\ell \ge b_\ell \ge c_\ell$ are the three radii along the three orthogonal axes and $d_\A$ is the estimated average diameter based on the observed magnitude, given in Table~\ref{t:apophis}. S-type asteroids, as used here for the test case, are moderately bright with an albedo from 0.10 to 0.22. By comparison, C-type asteroids are extremely dark with albedos typically in the range of 0.03 to 0.10. According to \citet{Delbo2007}, the geometric albedo of Apophis is 0.33 however the value used here of 0.2 was chosen to give a more general test case.

\begin{table}[htb] \begin{center}
   \caption{Orbital and physical properties of test asteroid.}   \label{t:apophis}
   \begin{tabular}{llll}
   \toprule
    \multicolumn{2}{c}{Element} &  Measured Value \\ \midrule
    Semi-major axis & $a_{\A_0}$ & 0.9224 AU \\
    Eccentricity & $e_{\A_0}$ & 0.1912 \\
    Inclination & $i_{\A_0}$ &  3.3312 deg \\
    RAAN & $\Omega_{\A_0}$ & 204.4428 deg \\
    Argument of periapsis & $\omega_{\A_0}$ &  126.4002 deg \\
    Period & $T_{\A_0}$ & 323.5969 days \\
    Mean motion & $n_{\A_0}$ & 1.2876 $\times 10^{-5}$ deg/s \\ \midrule
    Mass & $m_{\A_0}$ & 2.7$\times 10^{10}$ kg\\
    Gravitational constant & $\mu_\A$ &  1.801599$\times 10^{-9}$ km\spr{3}/s\spr{2} \\
    Physical dimensions & $a_\ell, b_\ell, c_\ell$ & 191 m,\;135 m,\;95 m\\
    Rotational velocity & $w_\A$ & 3.3$\times 10^{-3}$ deg/s \\
    Albedo & $\varsigma_\A$ & 0.2 \\  \bottomrule
\end{tabular}  \end{center} \end{table}

\begin{figure}[htb]
    \centering
    \parbox{0.48\textwidth}{ \includegraphics[width=0.5\textwidth]{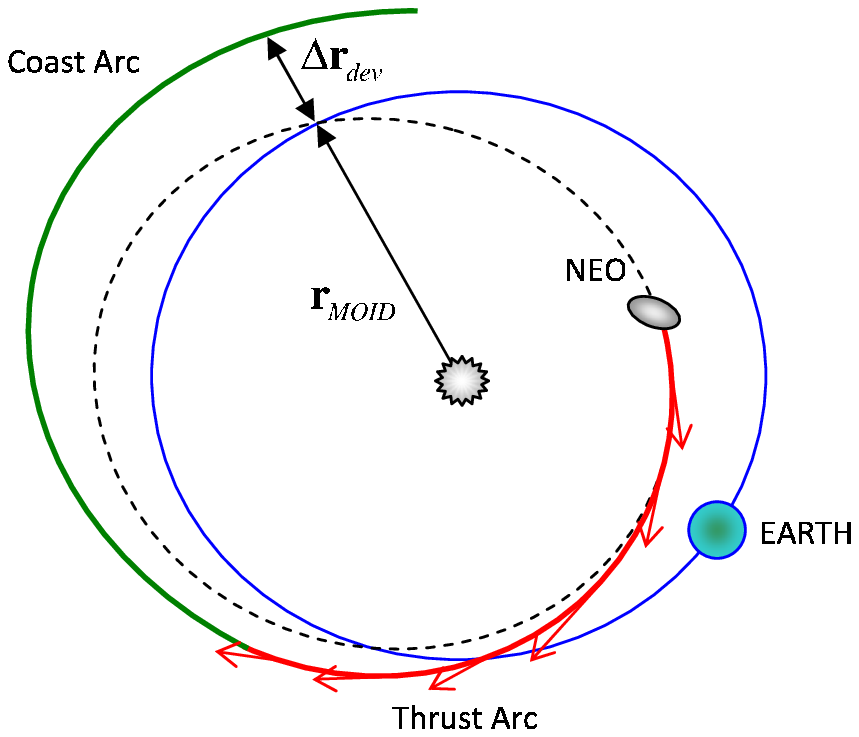}
    \caption{Definition of deviation distance at the MOID.}\label{fig:moid}}
    \quad
    \begin{minipage}{0.48\textwidth}
     \includegraphics[width=1\textwidth]{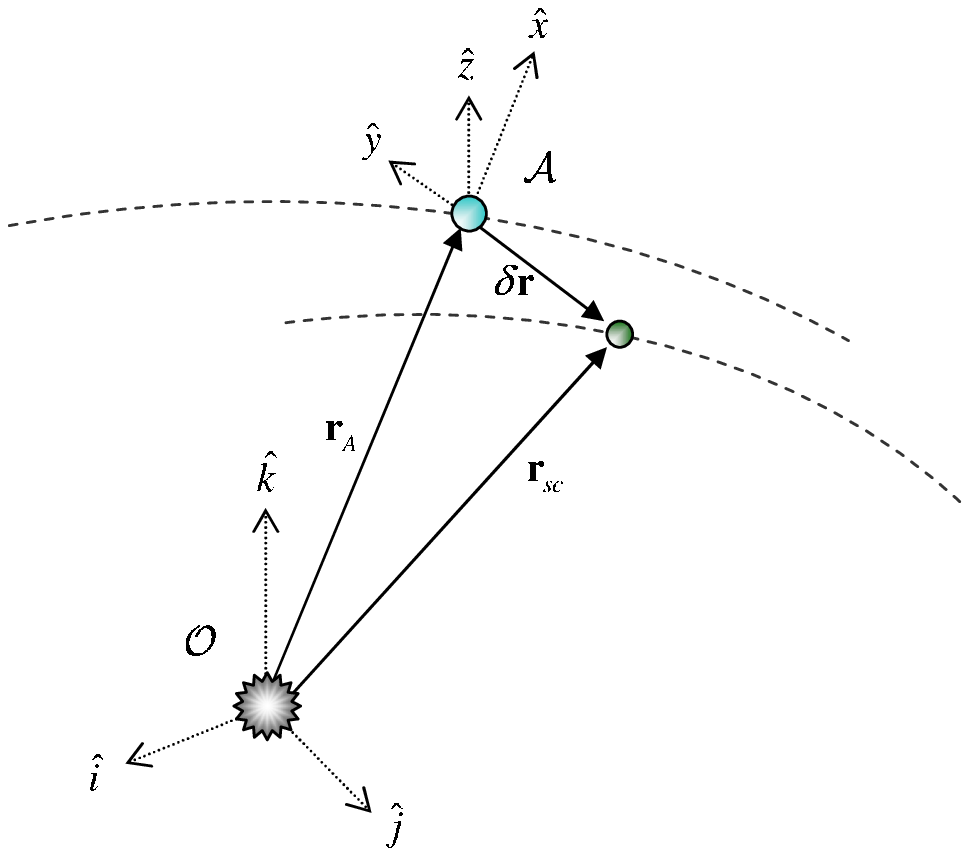} \caption{Definition of the reference frames, including the rotating Hill frame $\mathcal{A}$ centred on the asteroid.}\label{fig:hillframe}
    \end{minipage}
\end{figure}

The minimum orbital intersection distance (MOID) is the separation distance at the closest point between two orbits, e.g., Apophis and the Earth. The deviation distance is defined here as the difference in position between the original, undeviated orbit $\mathbf{k}_{\A_0}$ and the deviated orbit $\mathbf{k}_{\A_{dev}}$ at $t_\textsc{MOID}$ \citep{Colombo2009} (see Fig.~\ref{fig:moid}). Figure~\ref{fig:hillframe} illustrates the reference frames used here, where $\mathcal{O}\{i,j,k\}$ is the inertial heliocentric reference frame, and $\mathcal{A}\{x,y,z\}$ is the relative, rotating Hill reference frame (radial $x$, transverse $y$ and out-of-plane $z$ directions), centred on the asteroid.

Non-linear equations were used for determining the asteroid deviation vector $\Delta\mathbf{r}_{dev} = \mathbf{r}_{\A_{dev}} - \mathbf{r}_{\A_0}$ as a function of the ephemeris in the Hill reference frame $\mathcal{A}$, as derived by \citet{Maddock2008}, where $\Delta\mathbf{k}=\mathbf{k}_{\A_{dev}} - \mathbf{k}_{\A_0} = [\Delta a,\Delta e, \Delta i, \Delta\Omega, \Delta\omega, \Delta M]^\T$ giving the difference in Keplerian parameters between the undeviated and deviated orbits.

The change in the orbital parameters is calculated by numerically integrating the Gauss planetary equations \citep[see e.g.,][]{Battin} using a thrust vector $\mathbf{u}_\textit{dev}=[u_t\; u_n\; u_h]^\T$ in the tangential, normal and out-of-plane (or direction of angular momentum $h$) reference frame, induced by the deflection method:
\begin{equation}
\Delta\mathbf{k} = \int_{t_0}^{t_{MOID}} \frac{d\mathbf{k}(\mathbf{u}_\textit{dev})}{dt}\:dt
\end{equation}
Within this study, the deflection action is assumed to be aligned with the heliocentric velocity of the asteroid, therefore $u_n=0$ and $u_h=0$. Other authors have studied the optimal direction of the deflection action in the case of laser ablation \citep{Park2010}, however, the main interest of this paper is in the system sizing in relation to the achievable deviation.

\citet{Colombo2009} determined that the change in angular location, in this case given by the mean anomaly $M$, calculated at the MOID is,
\begin{equation}\label{eq:deltaM}
\Delta M = \int_{t_0}^{t_i} \frac{dM}{dt}\,dt + n_{\A_0}\left(t_0-t_\textsc{moid}\right)+n_{\A_i}\left(t_\textsc{moid}-t_i\right)
\end{equation}
where $n_{A_0}$ is the mean motion of the undeflected asteroid, $n_{A_i}$ is the mean motion of the asteroid at the end of the deflection action, $t_0$ is the beginning of the deflection action and $t_i$ is the end of the deflection action.

The non-linear proximal motion equations in \citet{vasile2010} together with \eqref{eq:deltaM} and the Gauss planetary equations give the variation of the orbit of the asteroid at the time of the MOID. \citet{Vasile2008a} showed that an estimation of the minimum orbit interception distance can be computed by projecting the variation of the orbit at the expected impact time onto the $b$-plane of the Earth at the time of the MOID, i.e., computing the variation of the impact parameter $b$. Hence, in the test section the variation of the impact parameter will be used as a measure of the achievable deflection.

The thrust produced by the deflection method is computed assuming that the lasers are not pulsed but continuous wave and that the energy density is sufficient only to turn the matter into gas (vapour regime) but not to produce plasma \citep{Phipps2010}. The level of momentum coupling that can be achieved with this model is lower than what can be found in other studies \citep[see e.g.,][]{Phipps2010}. A further assumption is that the asteroid is absorbing part of the incoming energy without changing its temperature thus providing a constant sink for heat transmission; this might not be the case for small asteroids.

Under these assumptions, the rate of the expelled surface matter is defined as \citep{Sanchez2009},
\begin{equation}\label{eq:dmdt}
\frac{dm_{exp}}{dt} = 2 n_{sc} v_{rot} \int_{y_0}^{y_{max}} \int_{t_{in}}^{t_{out}} \frac{1}{H} \left( P_{in}-Q_{rad}-Q_{cond}\right)  \; dt\; dy
\end{equation}
where $[t_{in}, t_{out}]$ is the duration for which a point is illuminated, $[y_0, y_{max}]$ are the vertical limits of the illuminated surface area (i.e.\ orthogonal to the direction of rotation of the asteroid), $H$ is the enthalpy of sublimation, $v_{rot}$ is the linear velocity of a point as it travels horizontally (i.e., orthogonal to $y$) through the illuminated spot area and $n_{sc}$ is the number of spacecraft in the formation.

The input power per unit area due to the solar concentrators is given by,
\begin{equation}\label{eq:pauPin}
 P_{in} = \eta_{sys} C_r (1-\varsigma_\A) S_0 \left( \frac{r_\AU}{r_\A} \right)^2
\end{equation}
where $\varsigma_\A=0.2$ is the albedo, $S_0=1367$~W/m\spr{2} is the solar flux at 1 AU, scaled to the Sun-asteroid distance $r_\A$, $\eta_{sys}$ is the system efficiency, and $C_r$ is the concentration ratio (the ratio between the power density from the Sun on the mirror surface, and that of the illuminated spot area on the asteroid).

The heat loss due to black-body radiation and the conduction loss are defined, respectively, as,
\begin{align}
Q_{rad} &= \sigma \epsilon_{bb} T^4   \label{eq:Qrad}\\
Q_{cond} &= (T_{subl}-T_0)\sqrt{\frac{c_\A k_\A \rho_\A}{\pi t}}\label{eq:Qcond}
\end{align}
where $\sigma$ is the Stefan-Boltzmann constant, $\epsilon_{bb}$ is the black body emissivity, $T$ is the temperature and $c_\A$, $\rho_\A$ and $k_\A$  are, respectively, the heat capacity, density and thermal conductivity of the asteroid. For the asteroid Apophis, $c_\A = 750$~J/kg$\cdot$K based on the average value for silicate materials, $k_\A=2$~W/K/m and $\rho_\A=2600$~kg/m\spr{3} \citep{Remo1994}. The sublimation temperature assumed is that for forsterites \citep{Wang1999}, $T_{subl} = 1800$ K, with $T_0$ set to 278~K.

The induced acceleration due to the sublimation process can then be determined by \citep{Sanchez2009},
\begin{equation}\label{eq:defact}
 \mathbf{u}_{sub} = \frac{\Lambda \,\overline{v}\,\dot{m}_{exp}}{m_\A} \, \mathbf{\hat{v}}_\A
\end{equation}
where $m_A$ is the mass of the asteroid at a generic instant of time, $\mathbf{\hat{v}}_\A$ is the direction of the velocity vector of the NEO, $\Lambda \simeq \left(\tfrac{2}{\pi}\right)$ is the
scattering factor, $\overline{v}$ is the average velocity of the debris particles according to Maxwell's distribution of an ideal gas:
\begin{equation}\label{eq:vexp}
    \overline{v}= \sqrt{\frac{8\mathrm{k}_\B T_{subl}}{\pi \textrm{M}_{Mg_2SiO_4}}}
\end{equation}
where $\mathrm{k}_\B$ is the Boltzmann constant, and $\textrm{M}_{\textit{Mg}_2SiO_4}$ is the molecular mass of fosterite.

The scattering factor $\Lambda$ is computed as the average of all possible thrust directions assuming that the thrust can point randomly at any angle $\alpha_\T$ between 0 and $\pi$, therefore $\Lambda=\tfrac{1}{\pi}\int_0^{\pi}\cos\alpha_\T \, d\alpha_\T$ \citep{Sanchez2009}. Some preliminary experiments \citep{gibbings2011} demonstrate that the plume is progressively focusing inwards for rocky type of asteroids, while for highly porous asteroids the plume tends to remain unfocused; hence assuming an uniform distribution of the thrust pointing direction over an angle of 180$^\circ$ is a conservative choice.
The remaining mass of the asteroid $m_\A$ is calculated by numerically integrating \eqref{eq:dmdt}.

\subsection{Contamination Model}

The contamination of the mirror surfaces due to the debris plume is modeled based on the work by \citet{Kahle2006}.  Their study made a number of initial assumptions regarding the expansion of the plume and sublimation process. The first assumption holds that the sublimation process is comparable to the generation of tails in comets. The asteroid is assumed to contain a reservoir of material underneath the surface, with the gas expanding outwards through a throat into vacuum. Preliminary experimental results have shown that this assumption, as with others in this section, are potentially overly pessimistic and may not be valid for every type of asteroid. However, altering these assumptions does not change the fundamental results in this paper, therefore it was decided to remain consistent with the existing literature and defer any further analysis on the validity of these assumptions for future work.

The second assumption is that the plume expansion is similar to the expansion of gas of a rocket engine outside the nozzle. The density of the expelled gas $\rho_{exp}$ is computed analytically,
\begin{equation}\label{eq:rhoexp}
    \rho_{exp}(\delta r_{s/sc},\varphi) = j_c \,\frac{\dot{m}_{exp}}{\overline{v}\,A_{spot}} \left(\frac{d_{spot}}{2 \delta r_{s/sc} + d_{spot}}\right)^2 \left(\cos\Theta\right)^{2/(\kappa-1)}
\end{equation}
where $d_{spot}$ is the diameter of the spot area, $\delta r_{s/sc}$ is the distance from the spot on the surface of the asteroid and the spacecraft, and $\Theta = \pi\varphi/2\varphi_{max}$ where $\varphi$ is the angle between the spot-spacecraft vector and the $y$-axis of the Hill reference frame. The jet constant $j_c$ was set to 0.345, the maximum expansion angle $\varphi_{max} = 130.45^\circ$, and adiabatic index $\kappa = 1.4$ based on the values for diatomic particles \citep{Legge1982}.

Note that this density model is in contradiction with the assumption of a uniform scattering over a hemisphere and, in fact, suggests a much more focused plume. From ongoing experiments \citep{gibbings2011}, the plume appears to more closely match the density distribution given in \eqref{eq:rhoexp} rather than a uniform distribution; nevertheless, in the analysis in this paper the most conservative choice was selected for the scattering factor in order to account for possible unmodeled performance degradation components.

The position vector $\delta\mathbf{r}_{s/sc}$ from the spot to the spacecraft is defined as:
\begin{equation}\label{eq:spot_distance}
  \delta\mathbf{r}_{s/sc}=\left[ \begin{array}{c}
    x-r_\ell \sin w_\A t \cos(-w_\A t- \theta_{v_\A})+r_\ell \cos w_\A t \sin(-w_\A t-\theta_{v_\A}) \\
    y-r_\ell \cos w_\A t \cos(-w_\A t- \theta_{v_\A})-r_\ell \sin w_\A t \sin(-w_\A t-\theta_{v_\A}) \\
    z
  \end{array} \right]
\end{equation}
where the radius of the ellipse is given by,
\begin{equation}\label{eq_r_ell}
    r_\ell=\frac{a_\ell b_\ell}{\sqrt{\big(b_\ell \cos(-w_\A t- \theta_{v_\A})\big)^2+\big(a_\ell \sin(-w_\A t- \theta_{v_\A})\big)^2}}
\end{equation}
and, with reference to Fig.\ \ref{fig:hillframe}, the position of the spacecraft with respect to the centre of the asteroid is $\delta \mathbf{r}=[x,\; y,\; z]^\T$. We assume here that the asteroid is spinning around the $z$ axis with a rotational velocity $w_\A$. The direction of the velocity of the asteroid in the heliocentric reference frame projected onto the Hill reference frame $\mathcal{A}$ is $\theta_{v_\A}$. In other words, in order to have a deflection thrust aligned with the velocity of the asteroid, the spot is assumed to be at an elevation angle over the $y$-axis equal to $\theta_{v_\A}$.

The third assumption made is that all the particles impacting the surface of the mirror condense and stick to the surface. The exhaust velocity is constant, therefore the thrust depends only on the mass flow. A higher thrust results in a higher mass flow and thus in a faster contamination. This is a rather conservative assumption. The actual contamination level depends on the type of deposited material and the temperature of the optical surfaces. Following the approach used to compute the contamination of surfaces due to out-gassing, a view factor $\psi_\textit{vf}$ was added equal to the angle between the normal to the mirror and the incident flow of gas. The resulting variation of the thickness of the material condensing on the mirror can be computed by,
\begin{equation}\label{eq:hcnd}
    \frac{dh_{cnd}}{dt} = \frac{2\, \overline{v}\, \rho_{exp}}{\rho_{_{layer}}} \cos \psi_\textit{vf}
\end{equation}
The average debris velocity $\overline{v}$ is multiplied by a factor of 2 to account for the expansion of the gas in a vacuum. The layer density $\rho_{_{layer}}$ was set to 1 g/cm\spr{3}. The power density on the asteroid surface is decreased based on the contamination of the mirrors.

A degradation factor $\tau$ is applied to the power beamed to the asteroid surface, based on the Lambert-Beer-Bouguer law \citep{Kahle2006},
\begin{equation}\label{eq:pwrdegred}
    \tau = e^{-2\upsilon h_{cnd}}
\end{equation}
where $\upsilon = 10^4$/cm is the absorption coefficient for forsterite. Note that the values of $\upsilon$ and $\rho_{_{layer}}$ are based on the assumption that the deposited material is dense and absorbs the light over the whole spectrum. This is again a rather conservative assumption; experiments have shown that while it appears to be valid for some silicates such as forsterite, this assumption may not hold true for all materials. As mentioned previously, further experimentation and analysis are underway, and will be the topic of future publications.

\eqref{eq:hcnd} is numerically integrated, along with the Gauss equations, for the period of the mission.

\subsubsection{Tugging Effect}
The spacecraft will fly in formation with the asteroid at a distance $\delta r$, thus exerting a gravitational pull on it \citep{gong2009}. The tugging acceleration $\mathbf{u}_{tug}$ is given by:
\begin{equation}\label{eq:tug_force}
    \mathbf{u}_{tug}= -n_{sc}\frac{G m_{sc}}{\delta r^2} \, \delta\mathbf{\hat{r}}
\end{equation}
where $G$ is the universal gravity constant and $m_{sc}$ is the mass of a spacecraft. The sum of $\mathbf{u}_{tug}$ and $\mathbf{u}_{sub}$ forms the total deflection acceleration  $\mathbf{u}_{dev}$. The acceleration $\mathbf{u}_{dev}$ is used with Gauss planetary equations in order to determine the change in the NEO orbit.

\subsection{The Laser System}

Lasers work on the general premise of exciting electrons by stimulating them with the addition of photons (or quantum energy), which temporarily boost them up to a higher energy state. This stimulation continues until a population inversion exists, where there are more electrons at a higher energy state, e.g., $E_1$ than at the lower (or original) state, e.g. $E_0$. The release of photons when the electrons drop back to their original base state produce an emission that, generally, has the same spectral properties of the stimulating radiation, and is therefore highly coherent. 
The energy that is not released as part of the output emission, is instead released as heat. This means that the laser must be continually cooled, which in space means large radiators.

In this paper two general methods of powering the laser are considered and defined as: \textit{direct pumping}, where the energy is directly used to excite the laser, and \textit{indirect pumping}, where an intermediate step is used to first convert the energy, e.g., solar radiation, into electricity.

Indirect solar-pumped lasers convert the solar energy first into electricity, which is then used to power the laser. Photovoltaic cells are an obvious choice for space applications. The drawback, of course, is the addition of an electrical power generator meaning added mass, size and power requirements. Direct solar-pumped lasers, by comparison, do precisely what the name suggests: the laser is directly energised using solar radiation. Due to the mismatch between the wide-band emissions of the Sun with the narrow absorption bands of lasers, the loss of available solar power is currently rather high. For example, the overlap between a Nd:YAG (neodymium-doped yttrium aluminium garnet) crystal absorption spectrum and the solar radiation spectrum is around 0.14 \citep{Weksler1988}.

One option is to use high efficiency solar arrays in conjunction with a solid state laser. Solid state lasers pumped with electric power can currently reach 60\% efficiency. If the solar arrays have an efficiency of 30\%, then the system would have an overall efficiency of 18\%. If a pumped laser is used, then the focal point can be close to the primary mirror and a high concentration factor can be obtained with a relatively small mirror. For example, if the mirror has an area of 314 m\spr{2} (equivalent to a 20 m diameter circular mirror), then the collected power at 1 AU is 429.5 kW. The solar array plus laser system converts only 18\% of this power, therefore only 77.3 kW is beamed to the surface to the asteroid, while the rest needs to be dissipated.

In a paper presented in 1994, Landis discussed the use of a directly solar pumped laser based on semiconductor technology. According to Landis, the expected efficiency of directly pumped semiconductor laser would depend on the same efficiency losses of a solar cell, therefore Landis was expecting a lasing efficiency (output/input power ratio) of 35\%. Such an efficiency would be one order of magnitude higher than the best Nd:YAG laser system, which is expected to reach 6\% of overall efficiency.

Direct solar pumping would represent an interesting solution in terms of complexity of the overall system. In fact no cooling system for the photovoltaic conversion and no power transmission would be required. On the other hand the Technology Readiness Level (TRL) of both solar cells and semiconductor lasers is far higher than the one of a directly pumped laser and an indirectly pumped laser can be expected to be operational much sooner.

Recent electrically pumped semiconductor laser have proven over 73\% wall-plug efficiency \citep{crump2005, stickley2005, nLIGHT, peters2007} with a target efficiency of 80\%. Research on fibres coupled with clusters of diodes have demonstrated slope efficiencies of up to 83\% \citep{Jeong2003, Jeong2004}. A substantial increase in cells efficiency has also to be expected. In particular, in order to achieve a 35\% efficiency in direct pumping, semiconductor technology should allow the absorbtion of the solar spectrum over a wide range of frequencies. A high efficiency of a directly pumped laser is therefore expected to correspond to a high efficiency of solar cells. An increase of solar cell efficiency up to 50\% \citep{Luque2004} is reasonable, allowing an indirect pumping system to have a comparable efficiency to a 35\% direct pumping system.

In the following the assumption is that the overall system efficiency $\eta_{sys}$ is about 22.7\%, with a 45\% efficiency of the cells, a 90\% efficient reflectors, a 85\% efficiency of the power transmission and regulation line and a 66\% efficiency of the laser (given by the product of the target 80\% for the laser diode and the achieved 83\% slope efficiency of the fibres). A second option with a 60\% laser efficiency and 40\% cell efficiency is also considered.

\section{Formation Design}

One idea for the orbital design is to have the spacecraft flying in formation with the asteroid, orbiting in tandem around the Sun (see Fig.\ \ref{fig:laser_scheme}). The spacecraft have to maintain their relative position with respect to the asteroid in order to keep the required power density on the same spot on the surface of the asteroid (note that the surface of the asteroid is moving under the spot light of the laser). Therefore, the formation orbits have to be periodic and in close proximity with relatively low excursion in the relative distance from the asteroid. On the other hand the spacecraft should minimise any impingement with the plume of debris and gas coming from the sublimation of the surface material.

In order to design the desired formation orbits, one can start by considering the local Hill reference frame $\mathcal{A}\{x,y,z\}$ in Fig.\ \ref{fig:hillframe} and the associated linearised version of proximal motion equations \citep{Schaub} used in the calculation of the asteroid deviation vector:
\begin{subequations}\label{eq:lin_proxy_s}
\begin{align}
x(\nu)&=\frac{a_Ae_A\sin\nu}{\eta}\delta M-a_A\cos\nu \delta e\\
y(\nu)&=\frac{r_\A}{\eta^3}(1+e_A\cos\nu)^2\delta M+r_\A\delta\omega+\frac{r_\A\sin\nu}{\eta^2}(2+e_A\cos\nu)\delta e+r_\A\cos i_A\delta\Omega\\
z(\nu)&=r_\A(\sin\theta_A\delta i -\cos\theta_A \sin i_A \delta\Omega)
\end{align}
\end{subequations}
where $\eta=\sqrt{1-e_A^2}$, $\theta_A=\nu+\omega_A$, $\nu$ is the true anomaly, $a_A,e_A,i_A,\omega_A$ are respectively the semi-major axis, eccentricity, inclination and argument of the perihelion of the orbit of the asteroid at a generic moment in time and $\delta \mathbf{k}=[\delta a,\delta e, \delta i, \delta\Omega, \delta\omega, \delta M]^\T$ are the variations of the orbital elements, with the imposed conditions $\delta r \ll r_\A$, and $\delta a=0$ in order to have periodic motion. These equations are a first approximation of the motion of the spacecraft and do not take into account the gravity field of the asteroid or solar pressure.

If the optimal thrust direction that maximises the deviation is along the unperturbed velocity vector of the asteroid \citep{Colombo2009}, then the exhaust gases will flow along the direction of the velocity of the asteroid projected in the Hill reference frame. Therefore, the position vector in the radial, transversal and out-of-plane reference frame was projected onto the tangential, normal, out-of-plane reference frame to give $\delta \mathbf{r}_{tnh}=[x_{tnh}, y_{tnh}, z_{tnh}]^\T$. Then, the size of the formation orbits projected in the $x_{tnh}$-$z_{tnh}$ plane was maximised. All the requirements on the formation orbits can be formulated in mathematical terms as a multi-objective optimisation problem with two objective functions,
\begin{subequations}\label{eq:moo_problem1}
\begin{align}
\min_{\delta\textbf{k}\in D} \max_\nu J_1 &= \delta r \\
\min_{\delta\textbf{k}\in D} \max_\nu J_2 &= -\arctan \left( \frac{\sqrt{x_{tnh}^2+z_{tnh}^2}}{y_{tnh}} \right)
\end{align} \end{subequations}
subject to the inequality constraint,
\begin{equation}\label{eq:moo_con}
C_{ineq} = \min_\nu | y(\nu)| -y_\textsc{lim} > 0
\end{equation}
where $y_\LIM$ is a minimum distance along the $y$-axis, and $D$ is the search space for the solution vector $\delta\textbf{k}$. Table~\ref{tab:formation_delta} defines the boundaries imposed on $D$. The boundary values were obtained by progressively increasing each of the boundaries from 0 to the value in the table, looking at the value of the maximum distance from the asteroid. Larger boundaries produce solutions with a better (lower) performance index $J_2$ but a higher performance index $J_1$.

\begin{table}
 \caption{Boundaries on the formation orbital parameters}\label{tab:formation_delta}
\begin{tabular}{cccccc}
\toprule
    & $\delta e$ & $\delta i$ & $\delta \Omega$ & $\delta \omega$ & $\delta M$ \\
    & ($10^{-7}$) & ($10^{-7}$ rad) & ($10^{-7}$ rad) & ($10^{-7}$ rad) & ($10^{-7}$ rad)\\  \midrule
  Lower bound & $-0.01$ & $-0.1$ & $-0.9$ & $-1.5$ & $-0.1$\\
  Upper bound & 0 & 0.1 & 0.9 & 1.5 & 0.5\\
  \bottomrule
\end{tabular}  \end{table}

\begin{figure}
 \begin{center}
 \subfigure[\label{fig:objf_deltaia}]{
  \includegraphics[width=0.5\textwidth]{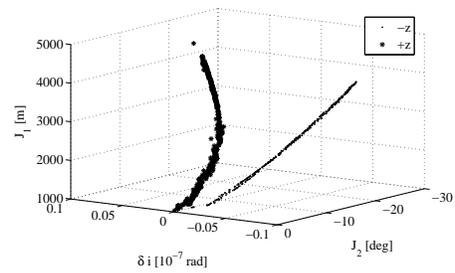}}
\subfigure[\label{fig:objf_deltaib}]{\includegraphics[width=0.5\textwidth]{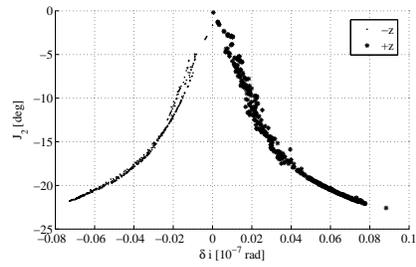}}\\
  \caption{a) Objective functions $J_1$ and $J_2$ versus $\delta i$, b) objective function $J_2$ versus $\delta i$}\label{fig:objf_deltai}
 \end{center}
\end{figure}

\begin{figure}
 \begin{center}
 \subfigure[\label{fig:objf_deltaOmegaa}]{
  \includegraphics[width=0.5\textwidth]{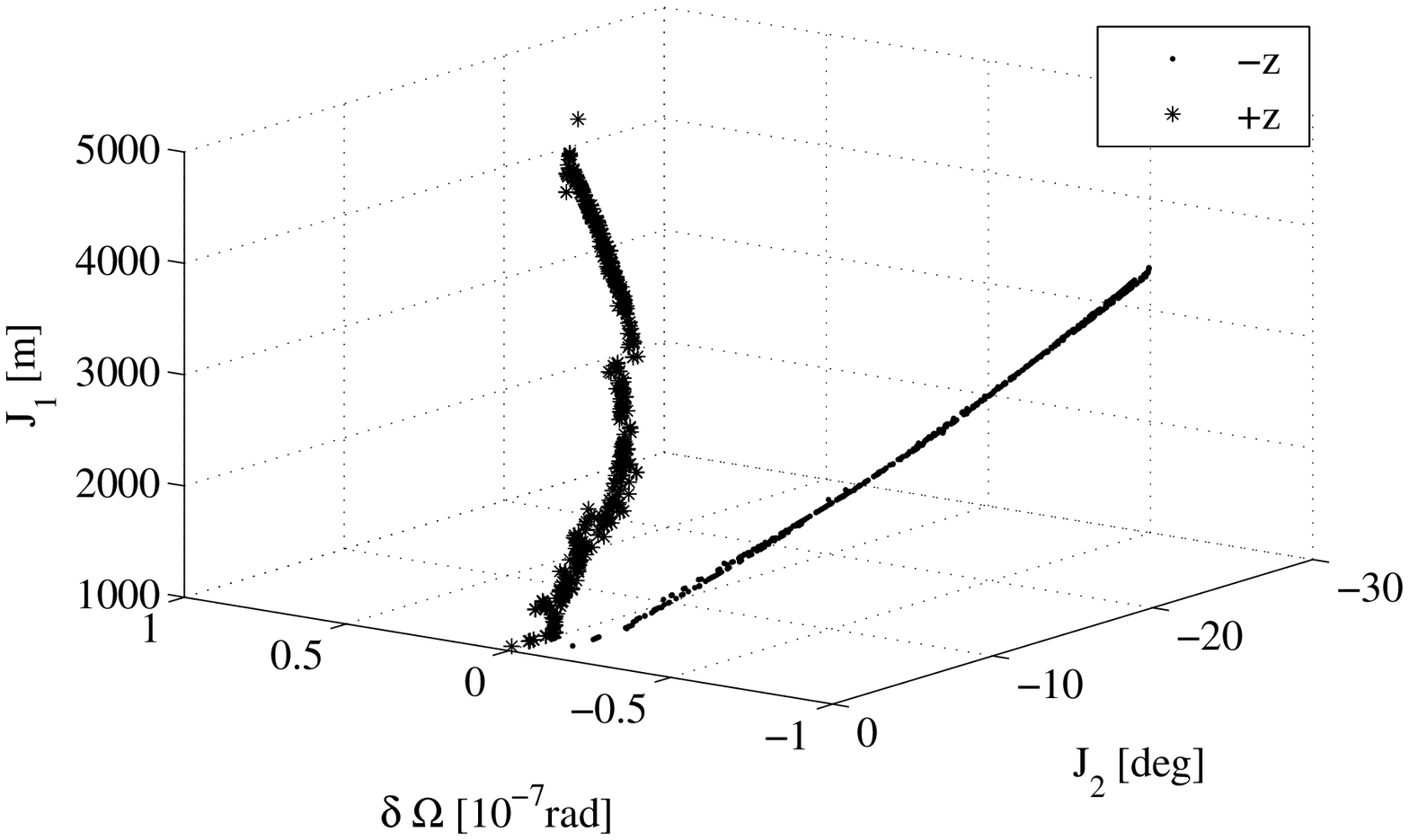}}
   \subfigure[\label{fig:objf_deltaOmegab}]{
  \includegraphics[width=0.5\textwidth]{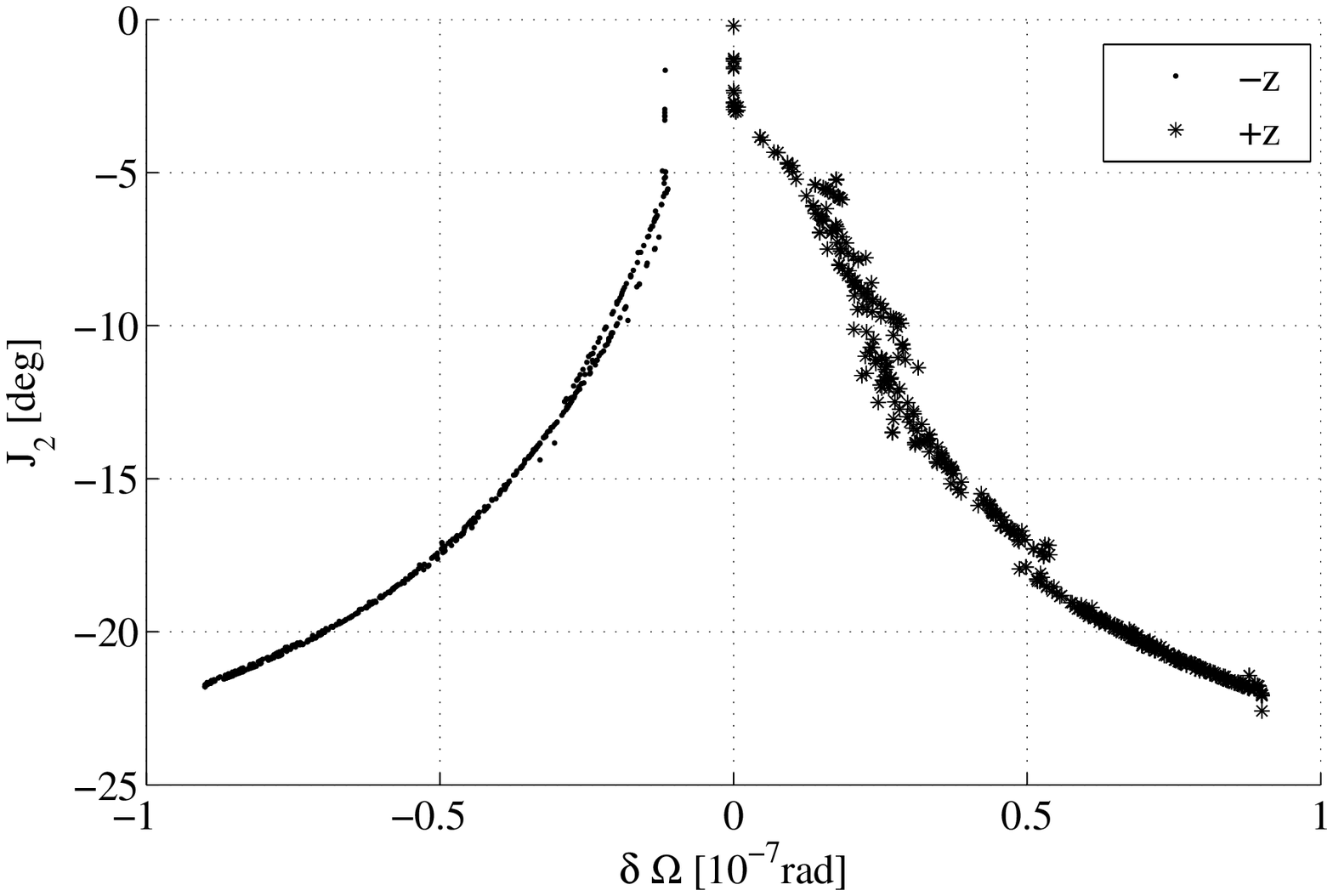}}\\
  \caption{a) Objective functions $J_1$ and $J_2$ versus $\delta \Omega$,b) objective function $J_2$ versus $\delta \Omega$}\label{fig:objf_deltaOmega}
 \end{center}
\end{figure}

\begin{figure}
 \begin{center}
 \subfigure[\label{fig:objf_deltaomegaa}]{
  \includegraphics[width=0.5\textwidth]{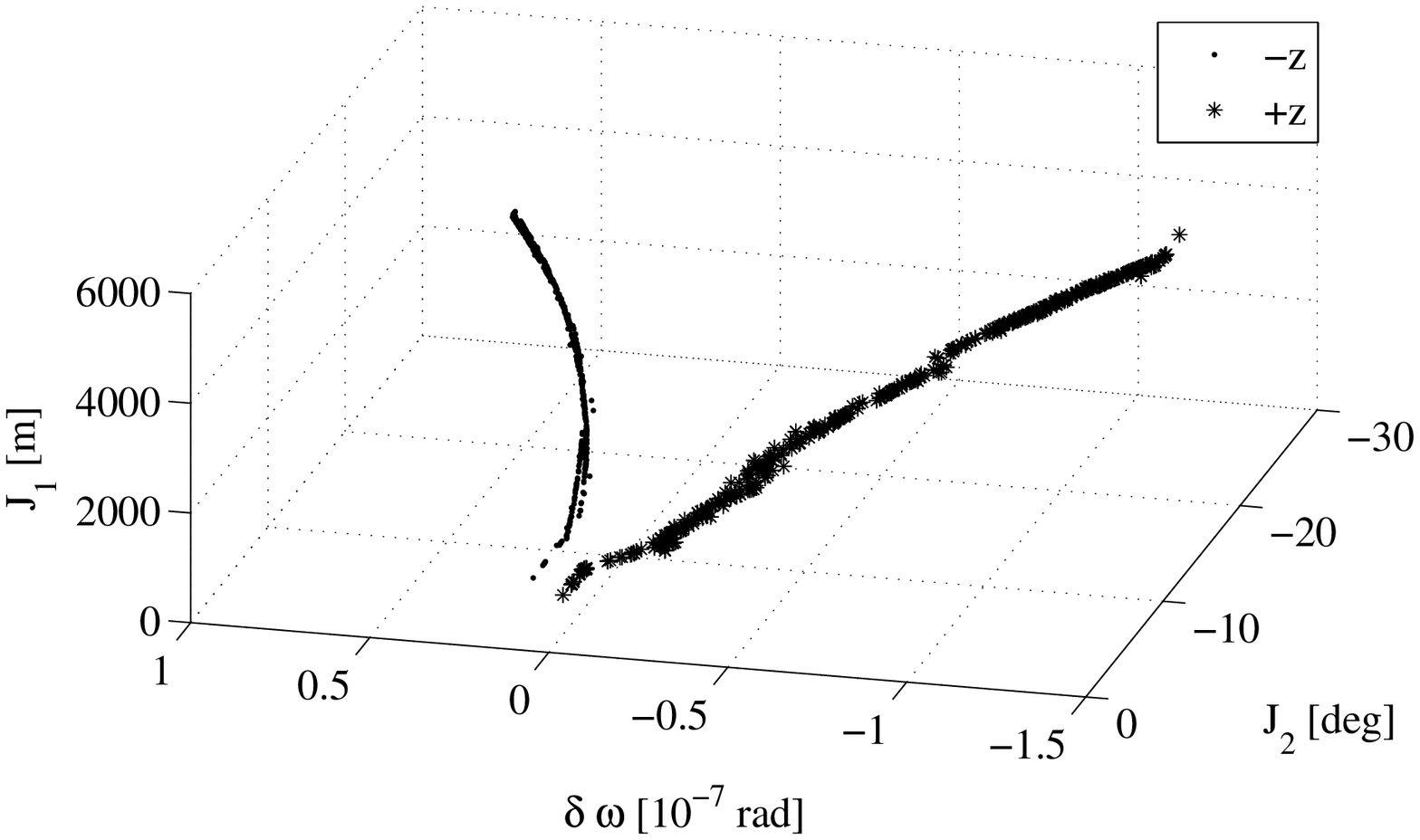}}
   \subfigure[\label{fig:objf_deltaomegab}]{
  \includegraphics[width=0.5\textwidth]{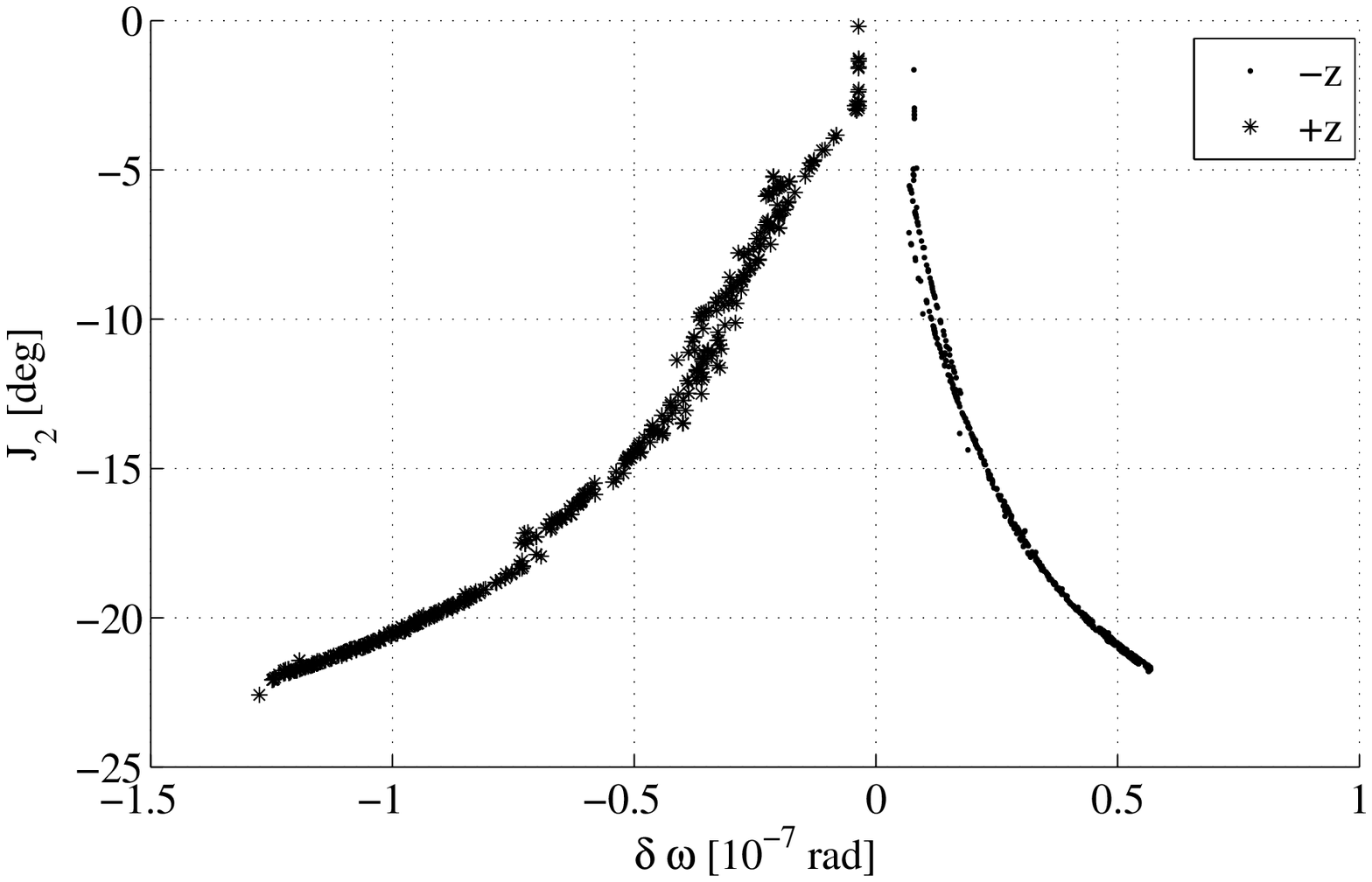}}\\
  \caption{a) Objective functions $J_1$ and $J_2$ versus $\delta \omega$, b) objective function $J_2$ versus $\delta \omega$}\label{fig:objf_deltaomega}
 \end{center}
\end{figure}

\begin{figure}
 \begin{center}
  \includegraphics[width=0.7\textwidth]{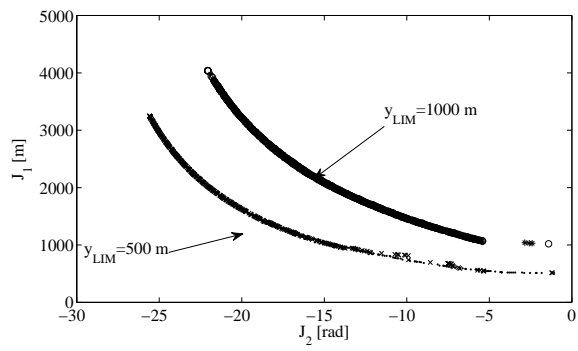}\\
  \caption{Pareto fronts of the optimal formation orbits}\label{fig:pareto_orbits}
 \end{center}
\end{figure}

Equations (\ref{eq:moo_problem1})--(\ref{eq:moo_con}) were optimised using a memetic multi-objective optimiser MACS (Multiagent Collaborative Search) \citep{Vasile2005a,Maddock2008,Vasile2008b}. The optimisation led to the identification of two families of formation orbits belonging to two subsets of the search space $D$. Figures \ref{fig:objf_deltai}, \ref{fig:objf_deltaOmega} and \ref{fig:objf_deltaomega} show the two families in the parameter space for $y_\LIM=1000$~m. The solutions are almost perfectly symmetrically distributed about the $0$-value of $\delta i$, $\delta \Omega$, while there is a bias towards the negative axis for $\delta \omega$.  Each family has been identified with the label $-z$ or $+z$ depending on whether the sign of the $z$ coordinate is negative or positive at $y=y_\LIM$. Figure~\ref{fig:pareto_orbits}, instead, shows the Pareto fronts for $y_\LIM=500$~m and $y_\LIM=1000$~m respectively. Note that in Fig.\ \ref{fig:pareto_orbits}, the Pareto fronts for the branches in Figs.~\ref{fig:objf_deltai}, \ref{fig:objf_deltaOmega} and \ref{fig:objf_deltaomega} appear superimposed and cannot be distinguished. Therefore, the two families can be considered equally locally Pareto optimal.

\begin{figure}[htb]
\begin{center}
  \includegraphics[width=0.7\textwidth]{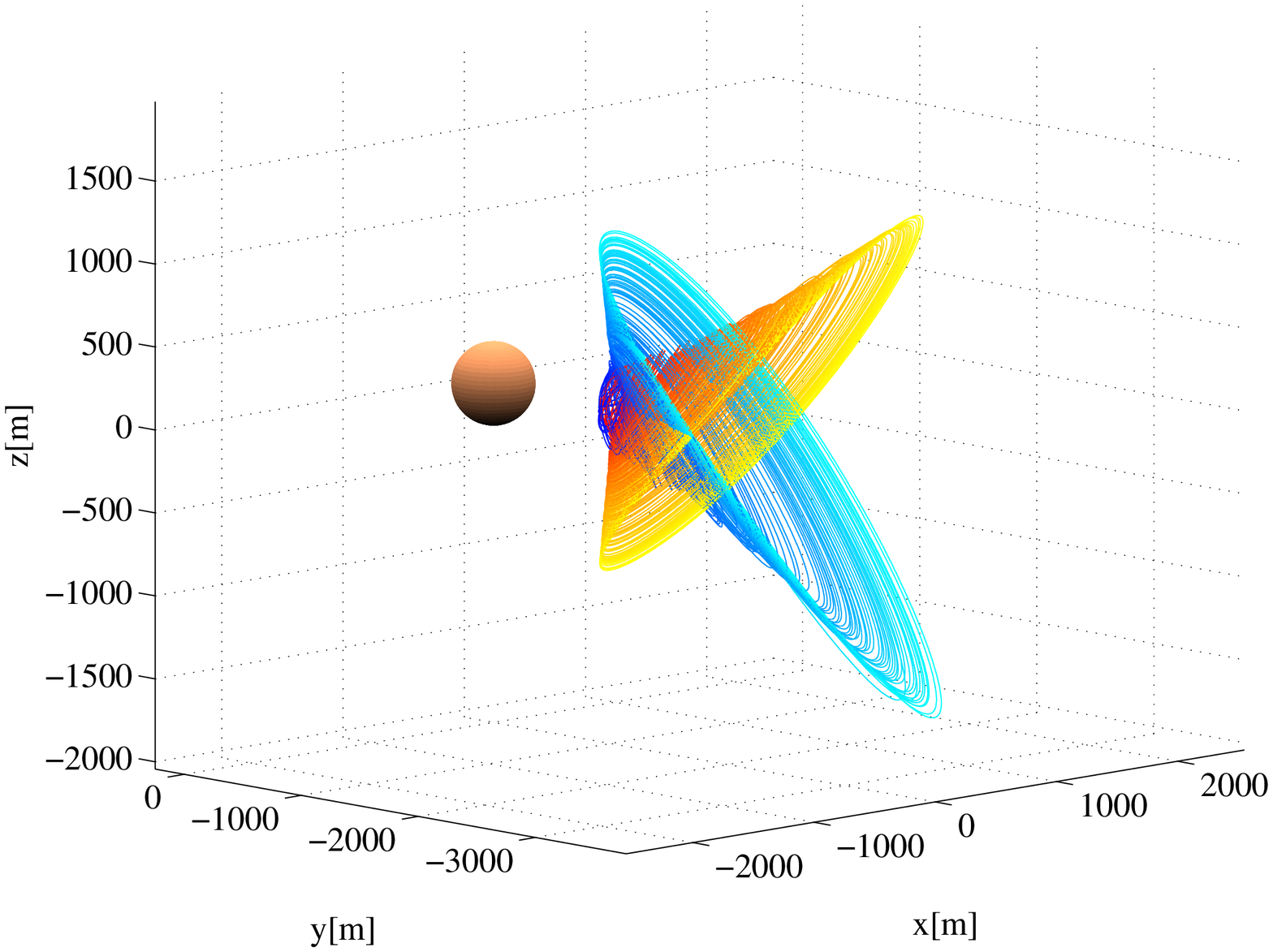}\\
  \caption{Formation orbits with minimum distance of 500 m.}\label{fig:form_orbit_min500}
\end{center}
\end{figure}

Figure \ref{fig:form_orbit_min500} shows the formation orbits in the $\mathcal{A}$ Hill frame. It can be noted that the two families are symmetric with respect to the $x$-$y$ plane. In the remainder of the paper these orbits will be called \emph{natural formation orbits}.

\subsection{Formation Dynamics and Control}

In order to maintain the orbits designed in the previous section, the spacecraft need to be controlled. In the proximity of the asteroid, in a Hill rotating reference frame, the spacecraft are subject to the force due to solar pressure, the gravity of the asteroid, the gravity of the Sun, the centrifugal and Coriolis forces plus the forces induced by the impingement with the plume. An active control is therefore required to maintain the spacecraft flying in formation with the asteroid.

Following the Jacobi ellipsoid model, the minor axis $c_\ell$ of the asteroid is aligned with vector of angular momentum, which corresponds to the $z$-axis of the asteroid Hill frame $\mathcal{A}$ (see Fig.\ \ref{fig:hillframe}).
The gravity field of the asteroid is expressed as the sum of a spherical field plus a second-degree and second-order field \citep{Hu2002,Rossi1999},
\begin{equation}\label{eq:grav_pot}
U_{20+22}=\frac{\mu_\A}{\delta r^3}\left(C_{20}\,(1-\frac{3}{2}\cos^2\gamma) + 3C_{22}\,\cos^2\gamma \cos 2\lambda\right)
\end{equation}
where $\gamma$ is the elevation over the $x-y$ plane and the harmonic coefficients $C_{20}$ and $C_{22}$ are a function of the semi-axes,
\begin{subequations}\label{eq:grav_coeff}
\begin{align}
C_{20}&=-\frac{1}{10}(2c_\ell^2-a_\ell^2-b_\ell^2)\\
C_{22}&=\frac{1}{20}(a_\ell^2-b_\ell^2)
\end{align}
\end{subequations}
and $\lambda$ is defined as,
\begin{equation*}\label{eq:sph_coor}
\lambda=\arctan\left(\frac{y}{x}\right)+ w_\A\, t
\end{equation*}
Note that a different rotational state or shape would alter the time-varying gravity field that the spacecraft would experience. In a real scenario, the rotational state coupled with the shape of the asteroid would require an adaptive focusing of the laser beam as the distance of the spot from the source will change with time. Also the divergence of the plume will change as the laser carves a groove into the asteroid. However, within the assumptions in this paper a different rotational state and/or shape would not alter the main results.

If one considers a Hill reference frame $\mathcal{A}$ centred in the barycentre of the asteroid (see Fig.\ \ref{fig:hillframe}), the motion of the spacecraft in the proximity of the asteroid is given by:
\begin{subequations}\label{eq:proxdynsolar3}
\begin{align}
\ddot{x}(\nu) &=-\ddot{r}_{A}+2\dot{\nu}\dot{y}+\dot{\nu}^2(r_\A+x)+\ddot{\nu}y - \frac{\mu_{\ssun}(r_\A+x)}{r_{sc}^3}- \frac{\mu_\A}{\delta r^3}x + \frac{F_{s_x}(x,y,z)}{m_{sc}} + \frac{\partial U_{20+22}}{\partial x}\\
\ddot{y}(\nu) &=-2\dot{\nu}\dot{x}-\ddot{\nu}(r_\A+x) +\dot{\nu}^2y - \frac{\mu_{\ssun}}{r_{sc}^3}y - \frac{\mu_\A}{\delta r^3}y + \frac{F_{s_y}(x,y,z)}{m_{sc}}+\frac{\partial U_{20+22}}{\partial y}\\
\ddot{z}(\nu) &=-\frac{\mu_{\ssun}}{r_{sc}^3}z - \frac{\mu_\A}{\delta r^3}z + \frac{F_{s_z}(x,y,z)}{m_{sc}}+\frac{\partial U_{20+22}}{\partial z} \label{eq:proxdynsolar_c}
\end{align} \end{subequations}
with,
\begin{align}
    \ddot{\nu} &= \frac{u_{dev_y}-2\dot{r}_\A r_\A\dot{\nu}}{r_\A^2}\label{eq:nuddot}\\
    \ddot{r}_\A &= \dot{\nu}^2 r_\A-\frac{\mu_\ssun}{r_\A^2}+u_{dev_x}\label{eq:rAddot}
\end{align}
The force term $\mathbf{F}_s=[F_{s_x}\, F_{s_y}\, F_{s_z}]^\T$ is made of two contributions: light pressure from the emitted light from the laser $\mathbf{F}_{srp}$  and the force due to the flow of gas and debris coming from the asteroid $\mathbf{F}_{plume}$.

The force due to solar radiation $\mathbf{F}_{srp}$ is defined as,
\begin{equation}\label{eq:force_srp}
    \mathbf{F}_{srp} =2\eta_{sys} A_{\M_1} \frac{S_0}{c} \left (\frac{r_\AU}{r_{sc}}\right)^2 \cos^2\beta \, \; \hat{\mathbf{n}}_{steer}+(1-\eta_\M^2) A_{\M_1} \frac{S_0}{c} \left (\frac{r_\AU}{r_{sc}}\right)^2 \hat{\mathbf{x}}
\end{equation}
where $c$ is the speed of light and $A_{\M_1}$ is the cross section area of the primary mirror (see \ref{fig:laser_scheme}). The angle $\beta$ is the half angle between the normal to the steering mirror $\hat{\mathbf{n}}_{steer}$ and the Sun-mirror vector (which is approximated by setting it equal to the Sun-asteroid vector). The second term in \eqref{eq:force_srp} takes into account a non-perfect reflection of the primary and secondary mirror. The reflectivity of the two mirrors is here assumed to be $\eta_\M=0.90$. The assumption is that the energy dissipated by the radiators is emitted uniformly in every direction and does not contribute to any change in the linear momentum of the spacecraft. 


If the flow rate per unit area at distance $\delta r_{spot}$ is $(2\rho_{exp}(\delta r_{spot},\varphi) \overline{v})$ and all the particles stick to the surface of the mirror then the force $\mathbf{F}_{plume}$ is:
\begin{equation}\label{eq:debries_force}
    \mathbf{F}_{plume}= 4\rho_{exp}(\delta r_{spot},\varphi) \bar{v}^2 A_{eq} \cos\psi_\textit{vf} \;\; \hat{\delta\mathbf{r}}_{s/sc}
\end{equation}
The flow rate depends on the power density and therefore on the distance from the Sun. The part of the spacecraft exposed to the plume and to the reflected light changes along the orbit and is irregular. In order to simplify the calculations, the assumption adopted in this paper is that the total effect is equivalent to a flat surface with area $A_{eq}=A_{\M_1}$ and normal unit vector $\hat{\mathbf{n}}_{eq}$ such that the cross product  $\langle\hat{\mathbf{n}}_{eq},\hat{\delta\mathbf{r}}_{s/sc}\rangle=\cos\psi_\textit{vf}$.

Given these equations, the resultant of all the forces acting on the spacecraft is not zero and in particular the difference between gravity and $\mathbf{F}_s$ is a function of time. Therefore, an active control is required to maintain the position of the spacecraft with respect to the asteroid.

If one assumes that solar pressure, the gravity of the asteroid, and the force due to the plume impingement are the main source of perturbation of the proximity motion of the spacecraft and that any non-spherical terms in the gravity field expansion result in only a small (second order) additional perturbation, then one can build a simple control law based on the Lyapunov control function:
\begin{equation}\label{eq:VAEP}
    V = \frac{1}{2}\delta v^2 + \frac{1}{2} K \left(  \left(x-x_\textit{ref}\right)^2 + \left(y-y_\textit{ref}\right)^2 + \left(z-z_\textit{ref}\right)^2  \right)
\end{equation}
where $\delta\mathbf{r}_\textit{ref} = [x_\textit{ref}, y_\textit{ref}, z_\textit{ref}]^\T$ are the coordinates of a point along the nominal formation orbit (in the Hill frame~$\mathcal{A}$). The assumption here is that the motion along the reference formation orbit is much slower than the control action, which is valid as the period of the spacecraft orbit is equal to the period of the asteroid (just under 1 year). Therefore, the spacecraft targets a set of static points along the formation orbit. Now if there exist a control $\mathbf{u}$ such that $dV/dt < 0$ then one can maintain the mirror in the proximity of the reference point as the reference point moves along the reference formation orbit. A possible control is given by:
\begin{equation}\label{eq:uVAEP}
    \mathbf{u} = - \left( -\frac{\mu_\A}{\delta r^3} \delta\mathbf{r} + \frac{\mathbf{F}_\textit{srp}}{m_{sc}}+ \frac{\mathbf{F}_\textit{plume}}{m_{sc}}\right)  - K \left(\delta\mathbf{r} - \delta\mathbf{r}_\textit{ref} \right) - c_d \delta\textbf{v}
\end{equation}
The total derivative of the function $V$ is:
\begin{subequations}\label{eq:dVdt}\begin{align}
\frac{dV}{dt} &= \delta\mathbf{v}^T \delta\mathbf{\dot{v}} + K(\delta\mathbf{r}-\delta\mathbf{r}_\textit{ref})^T  \delta\mathbf{v}\\
     & = \delta\mathbf{v}^T \Bigg( -\frac{\mu_\A}{\delta r^3}\delta\mathbf{r} + \frac{\mathbf{F}_\textit{srp}}{m_{sc}} + \frac{\mathbf{F}_\textit{plume}}{m_{sc}} -
     \left( -\frac{\mu_\A}{\delta r^3} \delta\mathbf{r} + \frac{\mathbf{F}_\textit{srp}}{m_{sc}}+ \frac{\mathbf{F}_\textit{plume}}{m_{sc}}  \right) \\ \nonumber
         &\qquad\qquad - K\left(\delta\mathbf{r} - \delta\mathbf{r}_\textit{ref} \right) - c_d \delta\mathbf{v} \Bigg) + K (\delta\mathbf{r} - \delta\mathbf{r}_\textit{ref})^T \delta\mathbf{v}\\
    & = -c_d \delta\mathbf{v}^T \delta\mathbf{v} < 0
\end{align}\end{subequations}
where $\delta\mathbf{v}=[\dot{x},\dot{y},\dot{z}]^\T$ is the relative velocity of the spacecraft in the asteroid Hill reference frame~$\mathcal{A}$.

The control in \eqref{eq:uVAEP} can now be introduced into the full dynamic model in \eqref{eq:proxdynsolar3} to test the validity of the assumption that the light coming from the asteroid and aspherical gravity field are indeed small. The elastic coefficient $K$ for both cases was chosen to be $10^{-6}$ while the dissipative coefficient $c_d$ was set to $10^{-5}$.

\begin{figure}[htb]
\begin{center}
  \includegraphics[width=0.7\textwidth]{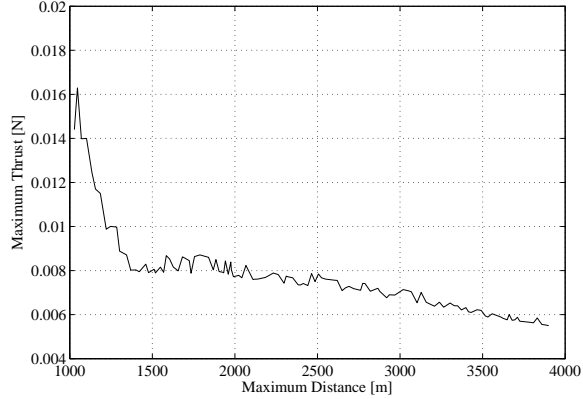}\\
  \caption{Maximum thrust versus maximum distance for the formation orbits with minimum distance of 1000 m.}\label{fig:maxT_maxD_min1000}
\end{center}
\end{figure}

\begin{figure}[htb]
\begin{center}
  \includegraphics[width=0.7\textwidth]{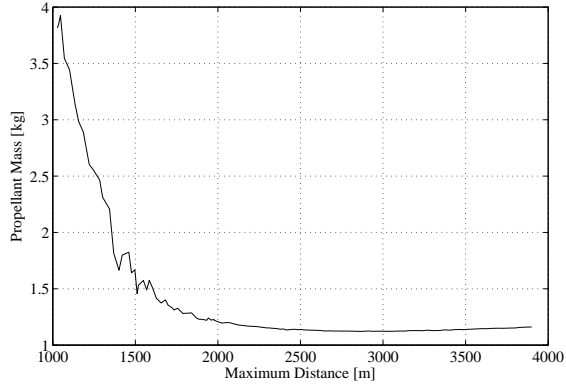}\\
  \caption{Propellant consumption versus maximum distance for the formation orbits with minimum distance of 1000 m.}\label{fig:prop_maxD_min1000}
\end{center}
\end{figure}

Figure \ref{fig:maxT_maxD_min1000} shows the maximum thrust level as a function of the maximum distance from the asteroid for a 20 m diameter mirror. Figure \ref{fig:prop_maxD_min1000} shows the propellant consumption as a function of the maximum distance from the asteroid for a 20 m diameter mirror.

\subsection{Shaped Formation}
Although the natural formation orbits are designed to minimise the impingement with the plume of gas and debris, none of them can avoid the plume completely. In order to maximise the amount of solar power collected, the mirrors should be constantly pointing directly towards the Sun, hence in a direction perpendicular to the $y$-axis. By following one of the natural formation orbits, the spacecraft will rise above the $z$-$y$ plane (i.e, in the $+x$ direction) once per revolution around the Sun, thus directly exposing the reflector to the plume. According to the contamination model, every surface directly exposed to the plume builds up a layer of a contaminants. This is quite a strong assumption as all the impinging material is assumed to condense and only the surfaces in view of the plume are contaminated. We hold on to these assumptions in this paper, although some experimental work is underway to build a more realistic model \citep{gibbings2011}. If one sticks to the assumptions of the contamination model, then one solution to mitigate the contamination would be to fly always below the plume of gas (i.e., $-x$ direction, below the $z$-$y$ plane).  In order to make the spacecraft follow the desired proximal motion the following shape is assigned to the formation orbit:
\begin{subequations}\label{eq:shaped_formations}
    \begin{align}
      x(\nu)&=x_1\cos(\nu)+x_2\sin(\nu)+x_3 \\
      y(\nu)&=y_1\cos(\nu)+y_2\sin(\nu)+y_3 \\
      z(\nu)&=z_1\cos(\nu)+z_2\sin(\nu)
    \end{align}
\end{subequations}
By differentiating with respect to time and inserting \eqref{eq:shaped_formations} and their first and second derivatives into the dynamic equations in \eqref{eq:proxdynsolar3}, one can compute the control profile and the corresponding propellant consumption. The interest now is to design formation orbits that minimise the propellant consumption required by the control system to remain below the $z$-$y$ plane and operate as close as possible to the asteroid to minimise pointing requirements. The problem can be formulated as follows:
\begin{subequations}\label{eq:moo_shaped_formation_problem}
\begin{align}
\min_{\textbf{s}\in X} J_1 &= \textit{MF}_\C\\
\min_{\textbf{s}\in X} \max_\nu J_2 &= \delta r\\
\min_{\textbf{s}\in X} \max_\nu J_3 &= \|\mathbf{u}\|
\end{align} \end{subequations}
subject to the inequality constraints:
\begin{subequations}\label{eq:moo_shaped_formation_con}
\begin{align}
C_{1} &= \max_\nu x(\nu)< 0\\
C_{2} &= \max_\nu y(\nu)< 0
\end{align}\end{subequations}
where the solution vector is $\mathbf{s} = [x_1, x_2, x_3, y_1, y_2, y_3, z_1, z_2]^\T$, and $\textit{MF}_\C$ is the propellant mass fraction for the control over one year of operations. The search space $X$ is defined by the lower and upper bounds on the components of $\mathbf{s}$, respectively $\mathbf{s}_{l}=[-1,    -1,   -1,    -1,  -1, -2, -1, -1]^\T $ and $\mathbf{s}_{u}=[1,     1,    0,     1,   1,  0, 1,  1]^\T $. Again MACS was used to solve the constrained problem in \eqref{eq:moo_shaped_formation_problem} and \eqref{eq:moo_shaped_formation_con}.

\begin{figure}[htb]
\begin{center}
  \includegraphics[width=1.1\textwidth]{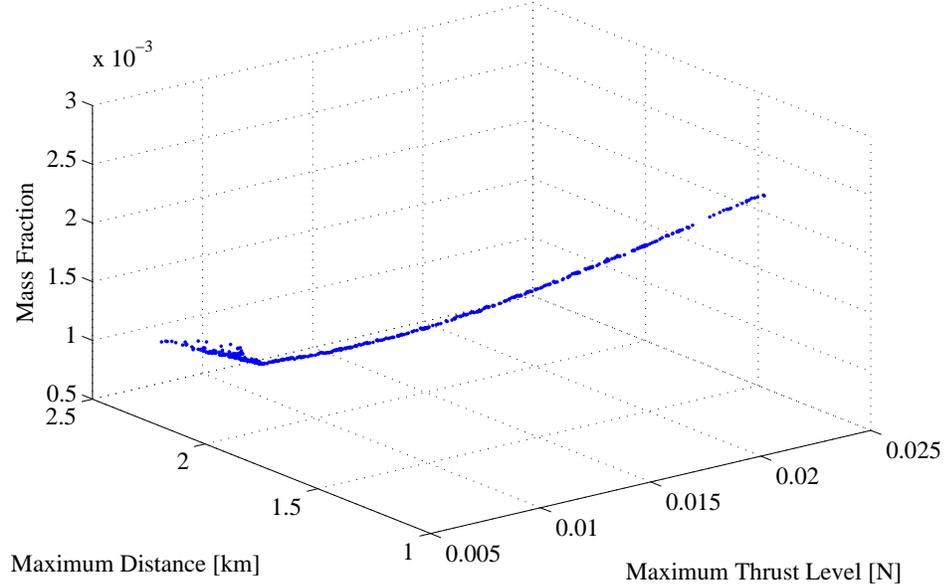}\\
  \caption{Pareto front of the shaped formation problem.}\label{fig:shaped_formation_pareto}
\end{center}
\end{figure}

The result of the multi-objective optimisation can be found in Fig.\ \ref{fig:shaped_formation_pareto}, and shows the propellant mass fraction versus the maximum thrust level versus the maximum distance to the asteroid for the case of 10 spacecraft, each carrying a 20~m diameter mirror, over the first year of operations.

As expected the level of thrust and control propellant mass fraction are monotonically increasing with the distance to the asteroid. However, even for close distances the annual propellant consumption and the thrust level are quite small, only a few milli-Newton of thrust are enough to maintain the orbit.

\section{Spacecraft and System Sizing}
\begin{figure}[htb]
  \centering  \includegraphics[width=0.9\textwidth]{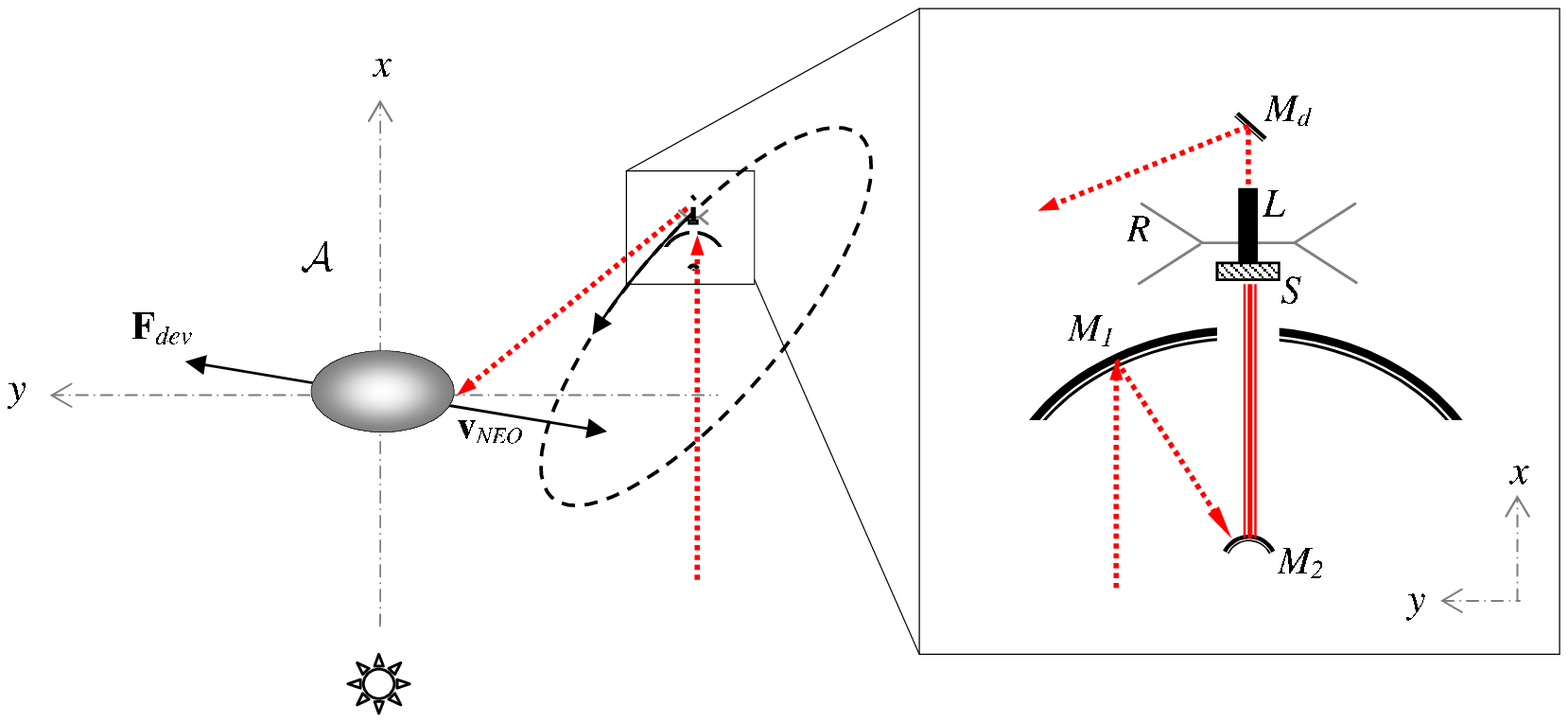}\\
  \caption{Illustration of the spacecraft and laser system, showing the two parabolic mirrors ($M_1$, $M_2$), the directional steering mirror ($M_d$), the solar arrays ($S$) which pump the laser ($L$), and the radiators ($R$).}
  \label{fig:laser_scheme}
\end{figure}

The proposed configuration of each spacecraft is as follows: each spacecraft is made of a primary mirror that focuses the sunlight onto a secondary mirror that reflects the light onto a solar array seated behind of the primary mirror (see Fig.\ \ref{fig:laser_scheme}). The electric power coming from the solar array pumps a semiconductor laser and a steering mirror directs the beam. The secondary mirror, the solar array and the laser need to be maintained at an acceptable temperature. Hence the need for radiators that dissipate the excess of energy that is not converted into the laser beam.

The size of the radiators can be computed considering the steady state thermal balance between the input power coming from the concentrator and the dissipated power through radiation.

Three radiating areas were considered for the design of the spacecraft: one associated to the secondary mirror with area $A_{R_{M2}}$, one associated to the solar array with area $A_{R_S}$, and one associated to the laser with area $A_{R_L}$. The size of each radiating area can be computed from the steady state equilibrium thermal equations:
\begin{subequations}\label{eq:radiatorA}\begin{align}
    A_{R_S}&=\left(\alpha_\s \eta_M P_{iM_2} -\eta_\s \eta_M P_{iM_2}-2\epsilon_\s \sigma A_{S} T_{S}^4\right) \Big/ \left(\epsilon_{\R} \sigma T_{\R_S}^4\right) \label{eq:S_radiator}\\
    A_{R_L}&=\eta_\s \eta_M P_{iM_2}(1-\eta_\La)\Big/ \left(\epsilon_\R \sigma T_\La^4\right) \label{eq:L_radiator}\\
    A_{R_{M2}}&=\left(\alpha_{\M_2}P_{iM_2}-2T_{\M_2}^4 \epsilon_\s \sigma \bar{A}_{\M_2}\right) \Big/ \left(\epsilon_\R \sigma T_{\M_2}^4\right)\label{eq:M_radiator}
\end{align}\end{subequations}
where $\alpha_\s$ is the absorptivity of the solar array, $A_S$ its area, $\epsilon_s$ its emissivity, $T_S$ its temperature, $P_{iM_2}=\eta_M A_{M_1} S_0 (r_{AU}/r_A)^2$ is the input power to the secondary mirror, and $\sigma$ is the Stefan-Boltzmann constant. Then, $\eta_S$ is the efficiency of the solar array, $T_{\R_S}$ is the temperature of the radiator associated to the solar array, and $\epsilon_{\R}$ its emissivity. Assuming the efficiency of the laser to be $\eta_\La$, and its temperature $T_\La$ one can compute the area of the radiator $A_{R_L}$ assuming that laser and radiator are in direct contact and that the heat can be transported with an efficiency close to 1. This is a reasonable assumption for relatively small scale systems that allow the use of a single or bi-phase passive cooling system. For large systems a dual phase active system might be required which lowers the overall efficiency and increases the system mass.
Finally, the secondary mirror is assumed to operate at temperature $T_{\M_2}$ and has absorptivity $\alpha_{\M_2}$.

\begin{table}\begin{center}
\caption{Thermal Properties of Spacecraft Elements}\label{tab:sc_elements}
\begin{tabular}{ccccccccc}
\toprule
\multicolumn{4}{l}{Solar arrays:} & \multicolumn{2}{l}{Mirror:} & Radiator: & \multicolumn{2}{l}{Laser:}\\
$\eta_\s$ & $\alpha_\s$ & $\epsilon_\s$ & $T_\s$ & $T_{\M_2}$ & $\alpha_{\M_2}$ & $\epsilon_\R$ & $\eta_\La$ & $T_\La$\\ \midrule
0.4-0.45 & 0.8 & 0.8 & 373 K  & 373 K& 0.01 & 0.9 & 0.6-0.66 & 313 K \\
\bottomrule
\end{tabular} \end{center}  \end{table}

\begin{table}\begin{center}
\caption{Mass of Spacecraft Elements}\label{tab:mass_sc_elements}
\begin{tabular}{lcc}
\toprule
  \hspace{4px}Specific mass: & $\varrho_{\M}$ &0.1 kg/m\spr{2}\hspace{10px}\\
    & $\varrho_{\M_d}$ &0.1 kg/m\spr{2} \\
    & $\varrho_\La$ &0.005 kg/W\\
    & $\varrho_\s$ &1 kg/m\spr{2}\\
    & $\varrho_\R$&1.4 kg/m\spr{2}\\    \midrule
  \hspace{4px}Mass: & $m_{bus}$& 500 kg\\ \midrule
  \hspace{4px}Mass fractions: & $\textit{MF}_\h$  & 0.2\\
     & $\textit{MF}_p$ &0.3 \\
     & $\textit{MF}_t$ &0.1 \\
    \bottomrule
\end{tabular} \end{center} \end{table}

The total mass of the spacecraft is $m_{sc}=m_{dry}+m_p(1 + \textit{MF}_t)$, where the mass of the propellant $m_p=m_{dry}\textit{MF}_p$ is a fraction of the dry mass $m_{dry}$, augmented by the mass fraction $\textit{MF}_t=10$\% to include the mass of the tanks.

The dry mass $m_{dry}=1.2(m_\h+m_\s+m_\M+m_\La+m_\R+m_{bus})$ is the sum of the mass of the laser $m_\La$, mass of the bus $m_{bus}$, mass of the mirrors $m_\M$, mass of the solar array $m_\s$, mass of the radiators $m_\R$ and mass of the harness $m_\h$. Given the low maturity of the technology employed for this system, we considered a system margin of 20\% on the dry mass.

The mass of the harness $m_\h$ is a fraction of the combined mass of the laser and solar array $m_\h=\textit{MF}_\h(m_\s+m_\La)$. The mass of the solar array is $m_s=1.15\varrho_\s A_\s$ where we considered a 15\% margin given the high efficiency of the cells.

The mass of the laser is $m_\La=1.5\varrho_\La P_\La \eta_\La$ where the margin is now 50\% given that a semiconductor laser of this size for space applications has not flown yet. The mass of the power management and distribution unit dedicated to the laser system is included in the mass of the harness while the mass of the bus is assumed to account also for the power electronics. The power input to the laser is,
\begin{equation}
P_\La=0.85\eta_\s \eta_\M^2 A_{M_1}S_0\left( \frac{r_\AU}{r_\A} \right)^2
\end{equation}
and is a function of the input light power on the solar array, the efficiency of the  solar array $\eta_\s$ and the reflectivity of the mirrors $\eta_\M=0.90$. The loss due to power regulation and transmission was considered to be 15\% of the generated power.

The mass of the radiators $m_\R=1.2(A_{\R_S}+A_{\R_{M_2}}+A_{\R_L})\varrho_\R$ from \eqref{eq:radiatorA} is proportional to the area and is augmented by a 20\% margin. The total mass of the mirror is $m_\M=1.25(\varrho_{\M_d} A_d + \varrho_{\M} \bar{A}_{\M_1} + \varrho_{\M} \bar{A}_{\M_2})$, where $\bar{A}_{\M_1}$ and $\bar{A}_{\M_2}$ are the areas of the primary and secondary mirrors. The total mass of the mirrors is augmented by a 25\% margin given the technology readiness level of the primary mirror.

The thermal properties of the system are reported in Table~\ref{tab:sc_elements} while the values of the specific masses $\varrho$, mass factors \textit{MF} and mass of the bus $m_{bus}$ are reported in Table~\ref{tab:mass_sc_elements}. The margins on the mirrors and power system are considered to include the marginal use of power to control the spacecraft in proximity of the asteroid. As it will be shown later, the required thrust level is small and therefore the power demand is marginal compared to the one required for the sublimation.

\subsection{Multiobjective Design}
Once the deflection and the spacecraft models are defined, the interest is to optimise the formation in order to obtain the maximum value of the impact parameter for the minimum mass into orbit, given a warning time. The problem can be formulated as follows:
\begin{subequations}\label{eq:moo_bees}
    \begin{align}
      &\min_{\mathbf{x}\in D} \left(-b\right)\\
      &\min_{\mathbf{x}\in D} \left(n_{sc}m_{sc}\right)
    \end{align}\end{subequations}
where the design vector $\mathbf{x}$ is defined by $[d_\M,\, n_{sc},\, C_r]^T$ and the search space $D$ is defined in Table~\ref{tab:D_space_mirror_design}.

\begin{table}\centering
  \caption{Boundary values on the design variables for the formation design.}\label{tab:D_space_mirror_design}
\begin{tabular}{lcc}
\toprule
    Design Parameter & Lower Bound & Upper Bound \\    \midrule
    Mirror aperture diameter, $d_\M$ (m) & 2 & 20 \\
    Number of spacecraft, $n_{sc}$ & 1 & 10 \\
    Concentration ratio, $C_r$ & 1000 & 5000 \\
    \bottomrule
\end{tabular} \end{table}

\begin{figure}[!htb]
\begin{center}
  \includegraphics[width=0.7\textwidth]{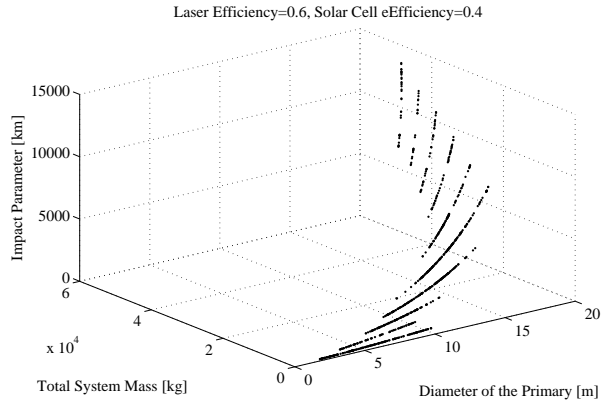}\\
\caption{Total mass of the system against the diameter of the primary mirror of each spacecraft and the achieved impact parameter. Natural formation orbits: $\eta_L=0.60$, $\eta_S=0.40$.}\label{fig:pareto_spacecraft_b_D_L6C4}
\end{center}
\end{figure}

\begin{figure}[!htb]
\begin{center}
  \includegraphics[width=0.7\textwidth]{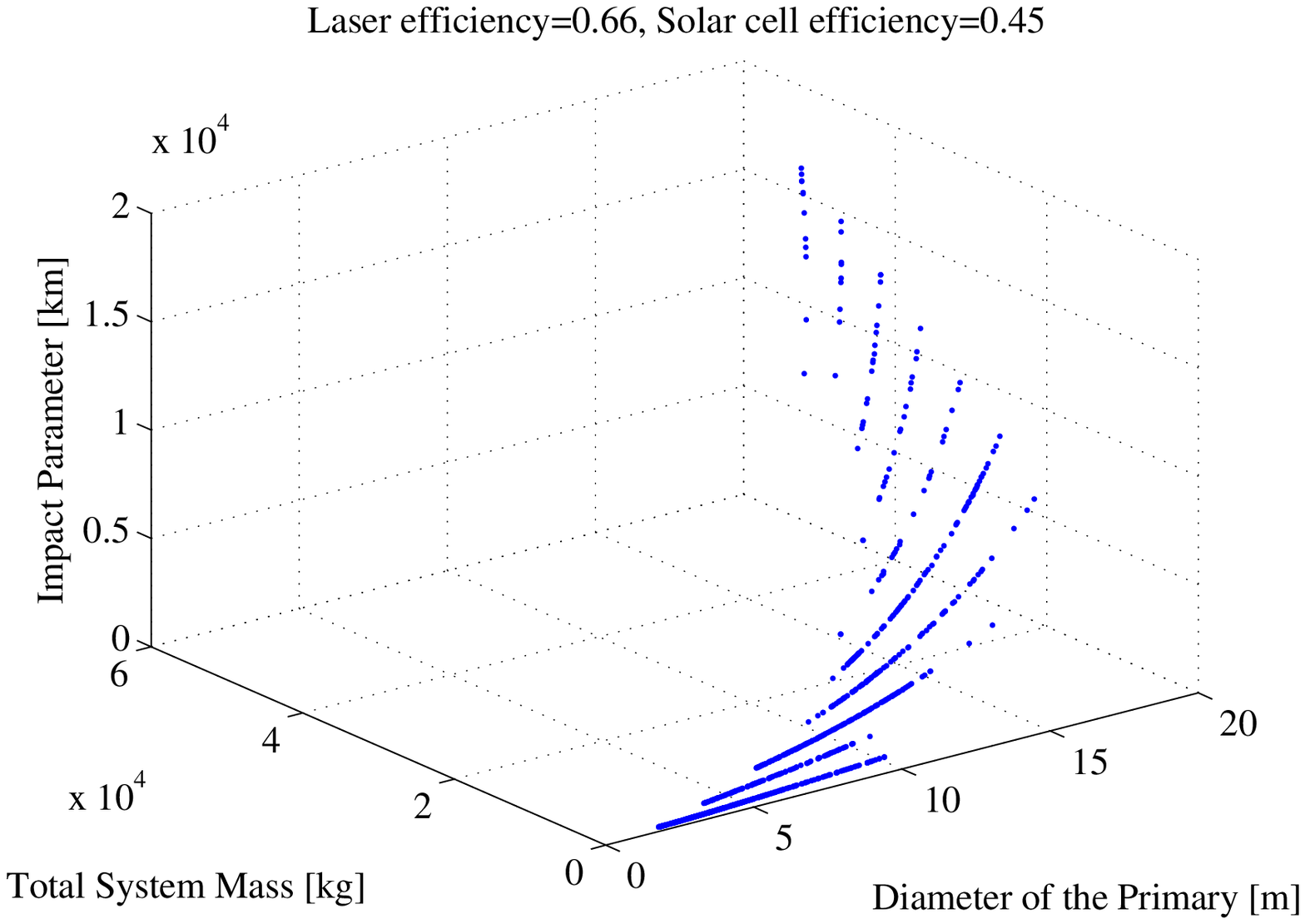}\\
  \caption{Total mass of the system against the diameter of the primary mirror of each spacecraft and the achieved impact parameter. Natural formation orbits: $\eta_L=0.66$, $\eta_S=0.45$.}\label{fig:pareto_spacecraft_b_D_L66C45}
\end{center}
\end{figure}

\begin{figure}[!htb]
\begin{center}
  \includegraphics[width=0.7\textwidth]{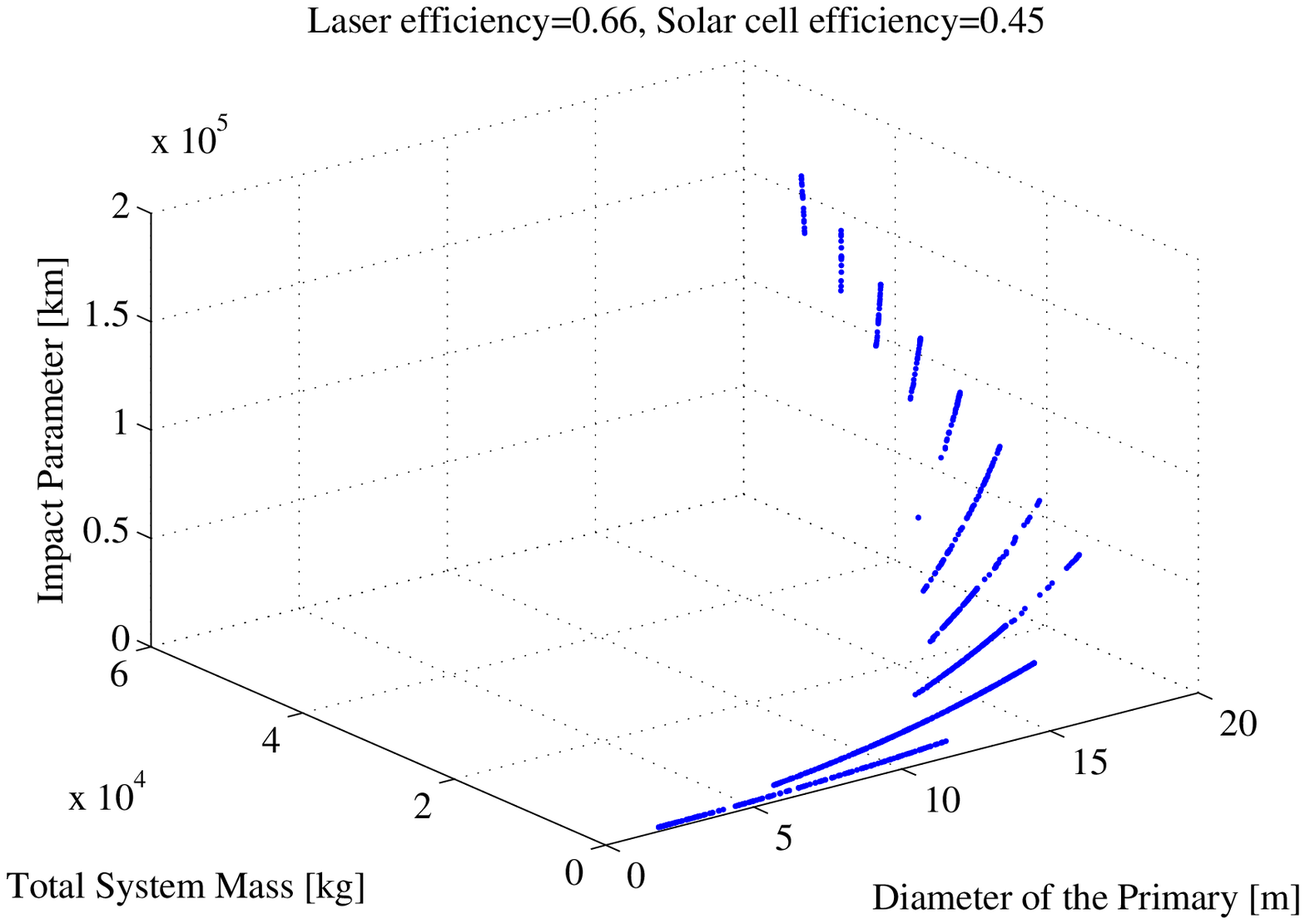}\\
  \caption{Total mass of the system against the diameter of the primary mirror of each spacecraft and the achieved impact parameter. Shaped formation orbits: $\eta_L=0.66$, $\eta_S=0.45$.} \label{fig:pareto_spacecraft_b_nsc_shaped_L66C45}
\end{center}
\end{figure}

\begin{figure}[!htb]
\begin{center}
  \includegraphics[width=0.7\textwidth]{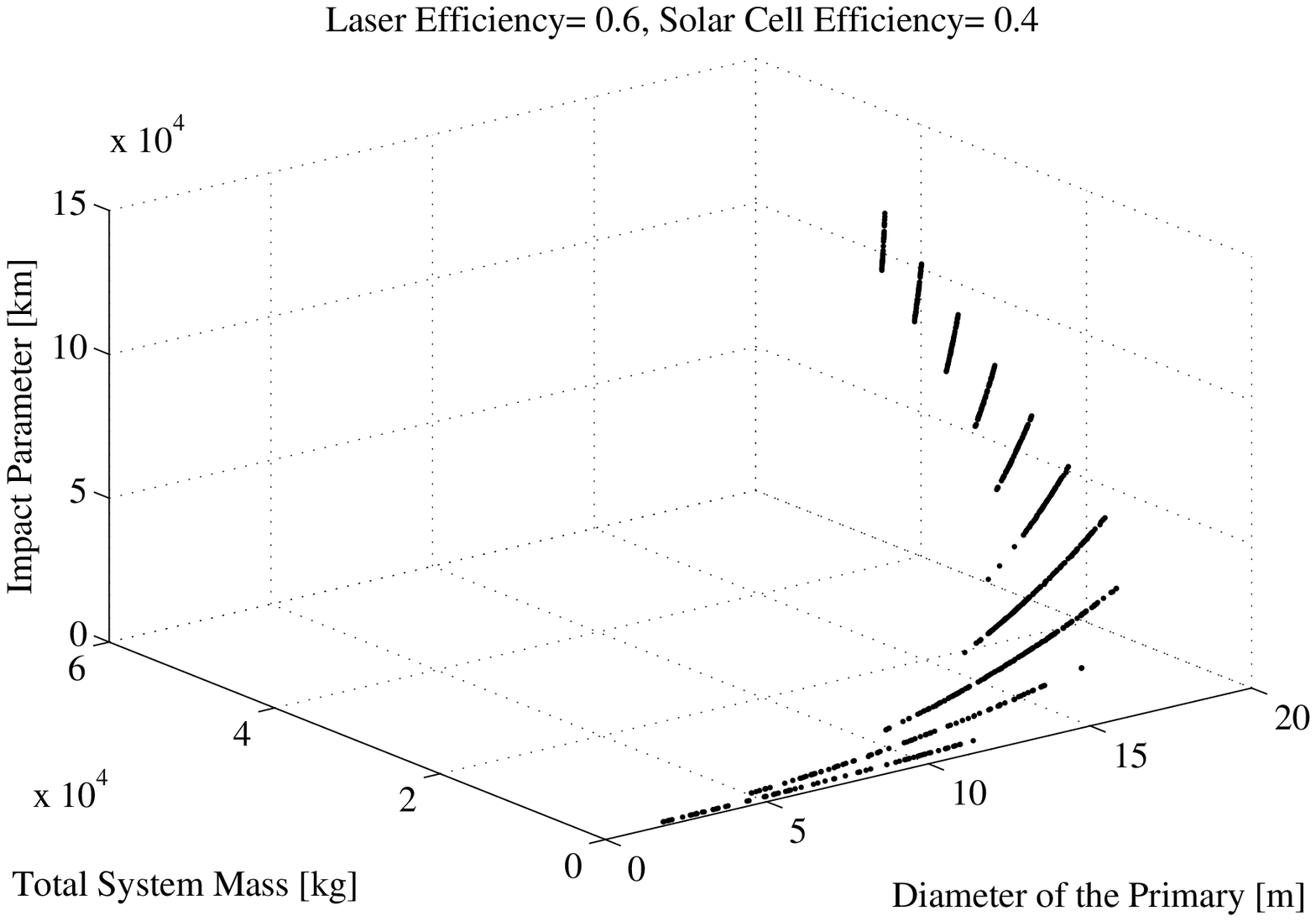}\\
  \caption{Total mass of the system against the diameter of the primary mirror of each spacecraft and the achieved impact parameter. Shaped formation orbits: $\eta_L=0.60$, $\eta_S=0.40$.} \label{fig:pareto_spacecraft_b_nsc_shaped_L6C4}
\end{center}
\end{figure}

The problem has two objectives, and a mix of integer and real variables. MACS was used to solve \eqref{eq:moo_bees}. The achievable deflection depends on the contamination of the optics, therefore the optimisation was run for both the shaped and the natural orbits. The result for the case of natural formation orbits can be seen in Fig.\ \ref{fig:pareto_spacecraft_b_D_L6C4}, where the impact parameter is represented against the mass of the system and the aperture diameter of the primary mirror for a laser with $\eta_L=0.6$ and solar cells with $\eta_S=0.4$, and Fig.\ \ref{fig:pareto_spacecraft_b_D_L66C45}, where the impact parameter is represented against the mass of the system and the aperture diameter of the primary mirror for a laser with $\eta_L=0.66$ and solar cells with $\eta_S=0.45$. Analogous solutions for the case of the shaped orbits can be found in Figs.\ \ref{fig:pareto_spacecraft_b_nsc_shaped_L66C45} and \ref{fig:pareto_spacecraft_b_nsc_shaped_L6C4}.

It is interesting to note that the number of spacecraft increases when the aperture diameter increases. This is due to that fact that as the diameter of the primary mirror increases the radiator and laser mass increases up to the point at which the mass of a single spacecraft exceeds the total mass of two or more spacecraft of smaller size. This is a very important point that is in favour of the use of a formation instead of a single large spacecraft.

Furthermore, it has to be noted that the assumption is that the system for each spacecraft is scalable. This is actually not true in general as the technology for radiators and concentrators cannot be arbitrarily scaled up. In other words, technological solutions for small size spacecraft cannot be applied to large size spacecraft without modifications. This is a further reason in favour of the use of multiple spacecraft of small size.

\begin{figure}[!htb]
\begin{center}
  \includegraphics[width=0.7\textwidth]{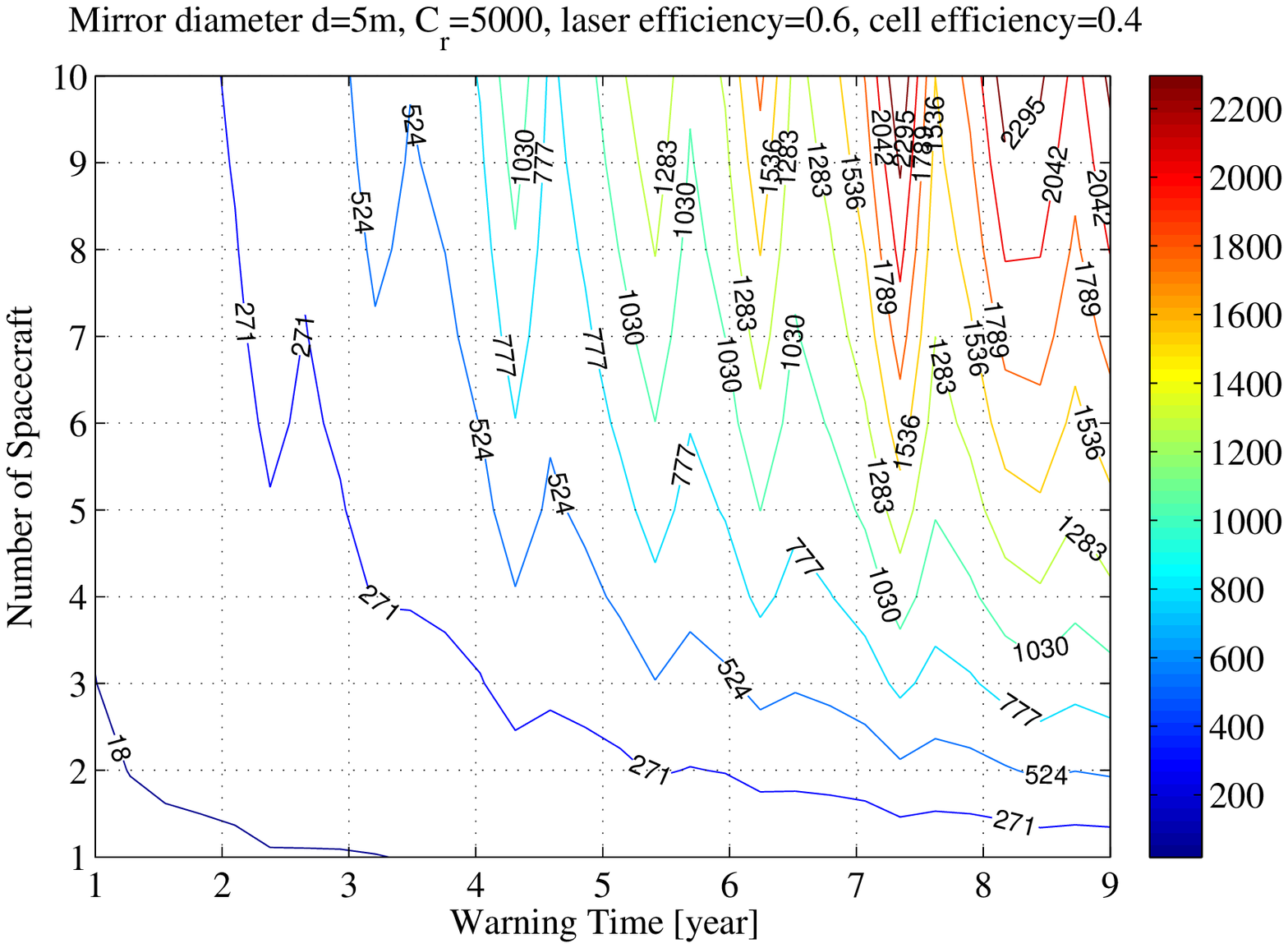}\\
  \caption{Impact parameter as a function of the number of spacecraft and warning time: 5 m aperture diameter and a concentration ratio of $C_r = 5000$. Natural formation orbits: $\eta_L=0.60$, $\eta_S=0.40$.} \label{fig:deflection_5mlaser}
\end{center}
\end{figure}

\begin{figure}[!htb]
\begin{center}
  \includegraphics[width=0.7\textwidth]{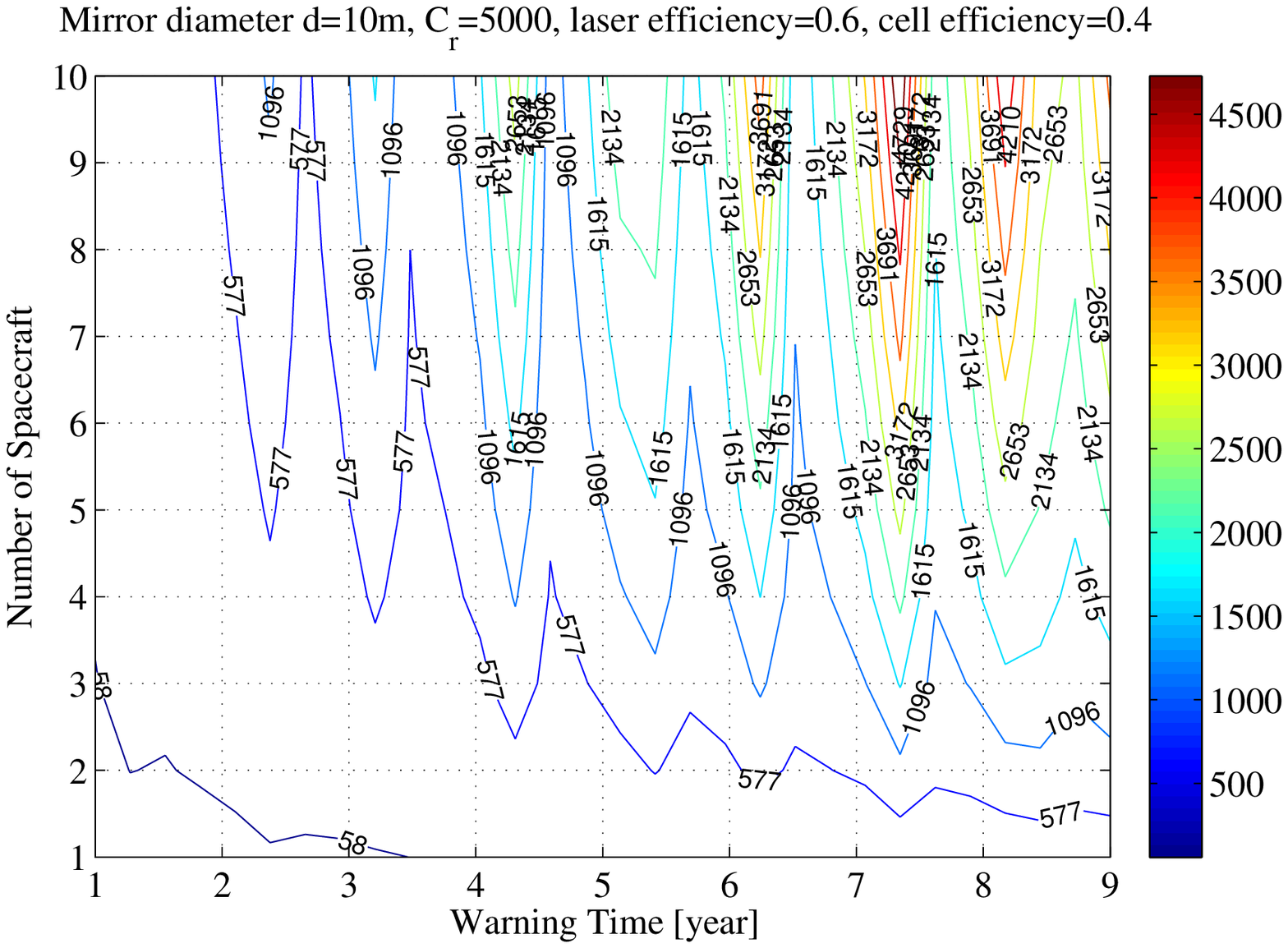}\\
  \caption{Impact parameter as a function of the number of spacecraft and warning time: 10 m aperture diameter and a concentration ratio of $C_r=5000$. Natural formation orbits: $\eta_L=0.60$, $\eta_S=0.40$.} \label{fig:deflection_10mlaser}
\end{center}
\end{figure}

\begin{figure}[!htb]
\begin{center}
  \includegraphics[width=0.7\textwidth]{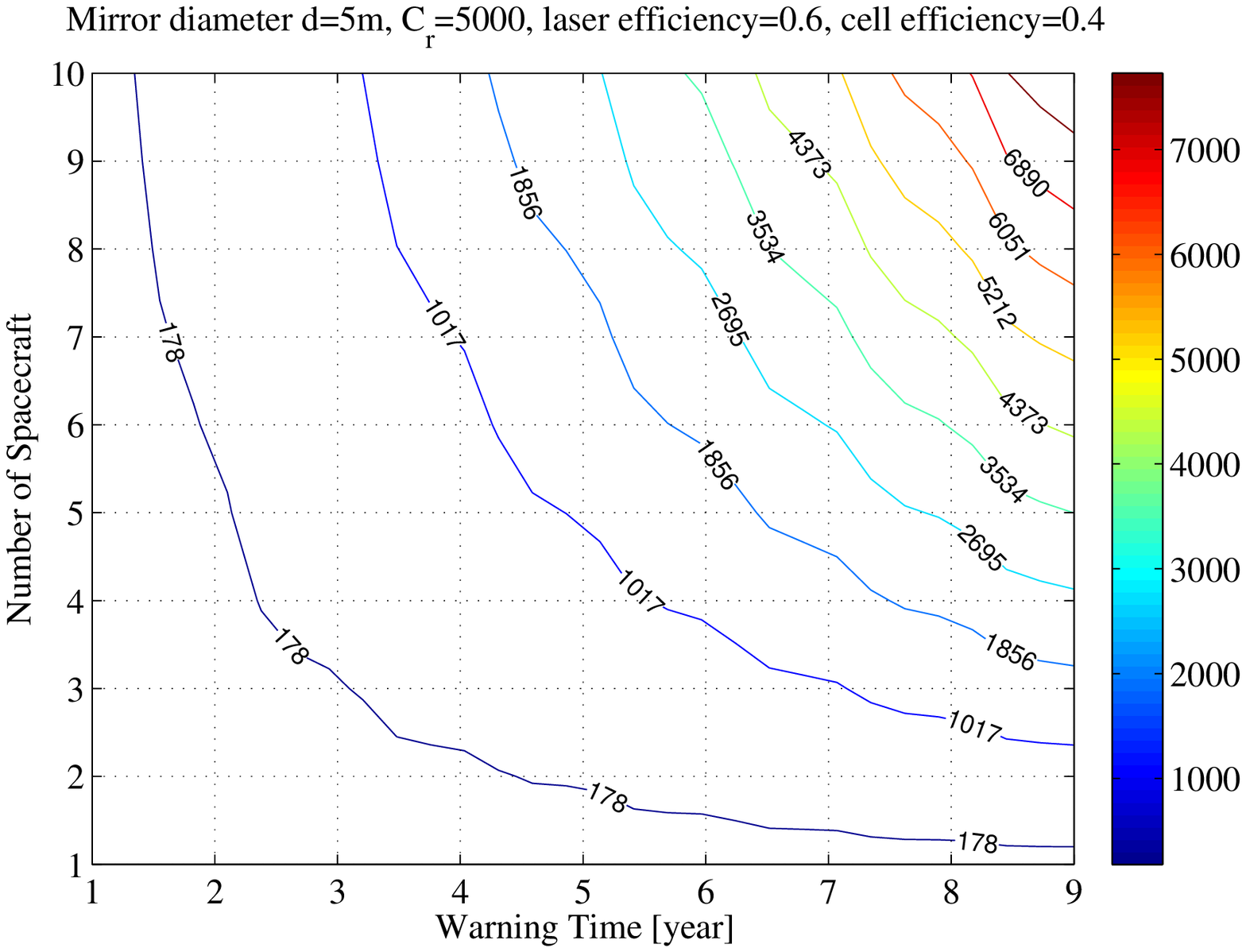}\\
  \caption{Impact parameter as a function of the number of spacecraft and warning time: 5 m aperture diameter and a concentration ratio of $C_r = 5000$. Shaped formation orbits: $\eta_L=0.60$, $\eta_S=0.40$.} \label{fig:deflection_shaped_5mlaser_shaped}
\end{center}
\end{figure}

\begin{figure}[!htb]
\begin{center}
  \includegraphics[width=0.7\textwidth]{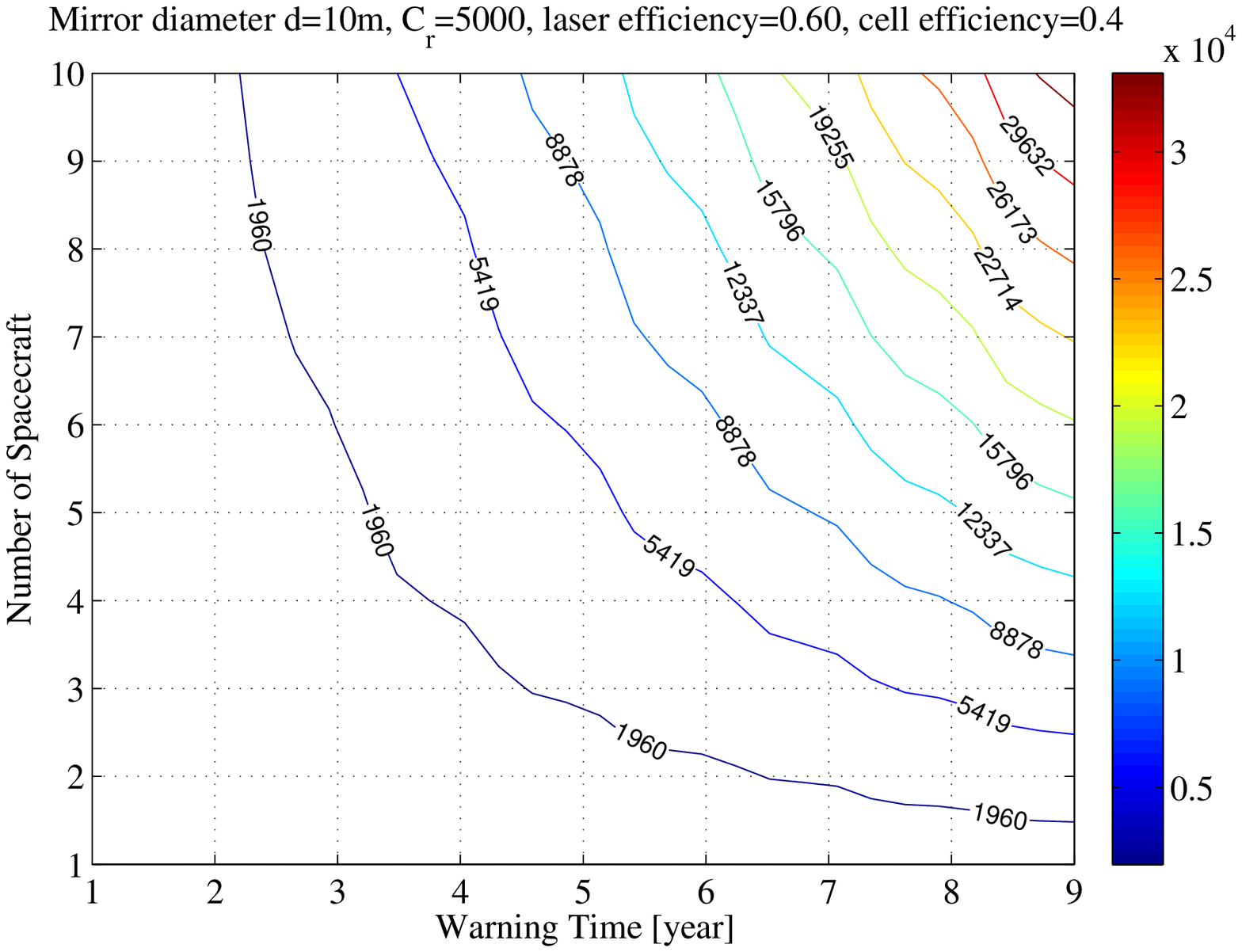}\\
  \caption{Impact parameter as a function of the number of spacecraft and warning time: 10 m aperture diameter and a concentration ratio of $C_r=5000$.  Shaped formation orbits: $\eta_L=0.60$, $\eta_S=0.40$.} \label{fig:deflection_shaped_10mlaser_shaped}
\end{center}
\end{figure}

\begin{figure}[!htb]
\begin{center}
  \includegraphics[width=0.7\textwidth]{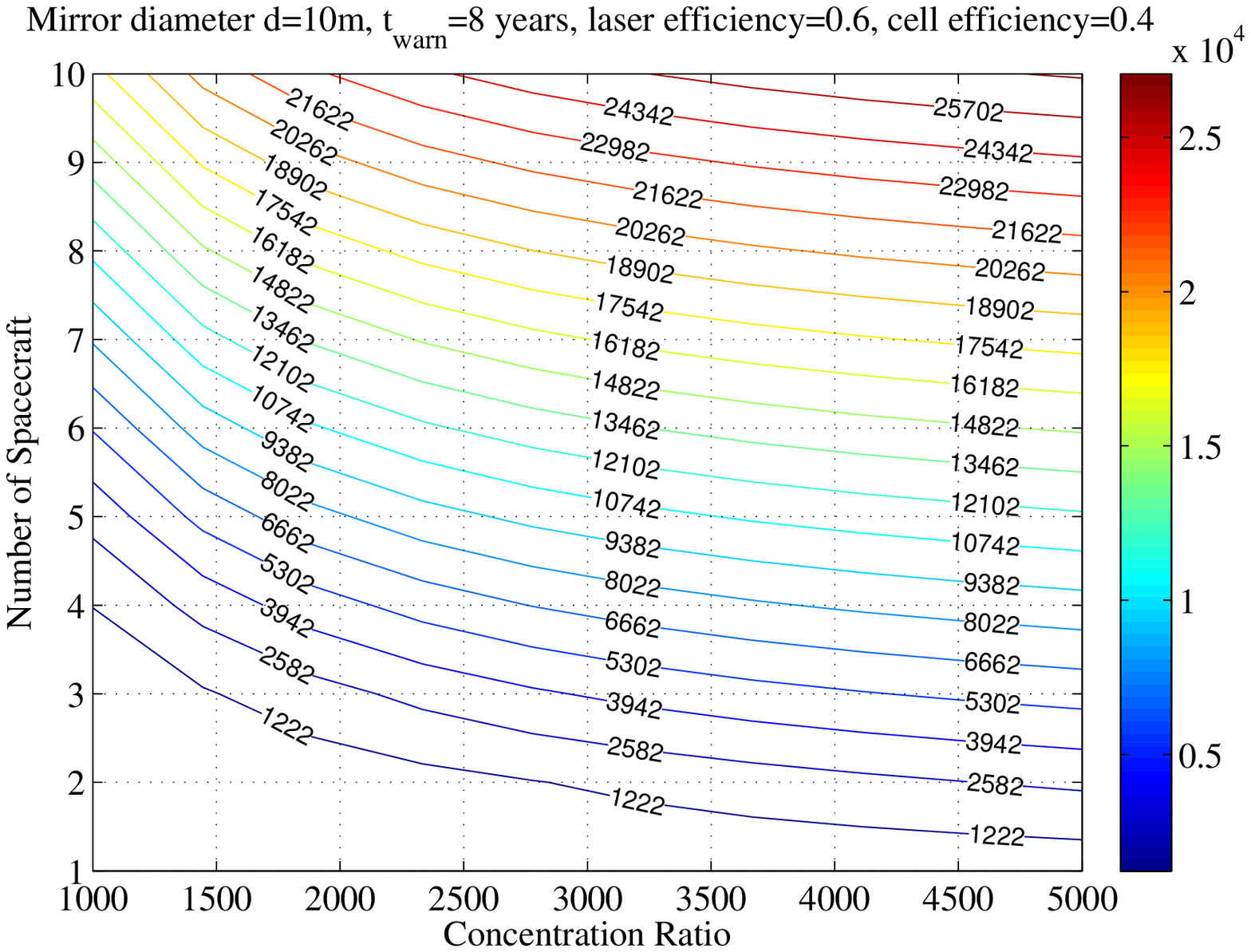}\\
  \caption{Impact parameter as a function of the number of spacecraft and the concentration ratio: 10 m aperture diameter and a warning time of 8 years.  Shaped formation orbits: $\eta_L=0.60$, $\eta_S=0.40$.} \label{fig:deflection_shaped_Cratio_shaped}
\end{center}
\end{figure}

Figures \ref{fig:deflection_5mlaser} and \ref{fig:deflection_10mlaser} show the achievable impact parameter for the case of the natural formation orbits with two alternative design solutions, a 5~m in diameter reflector and a 10m in diameter reflector both with a concentration ratio of 5000, i.e., the ratio between the area of the concentrator and the area of the spot on the surface of the asteroid is 5000. Figures \ref{fig:deflection_shaped_5mlaser_shaped} and \ref{fig:deflection_shaped_10mlaser_shaped} show the achievable impact parameter for the case of the shaped formation orbits. Figure \ref{fig:deflection_shaped_Cratio_shaped} shows the sensitivity to the concentration ratio for a fixed warning time of 8 years. The evident difference between the achievable impact parameter in the case of natural and shaped formation depends on the contamination effect that stops the sublimation process quite rapidly (less than one year in some cases) when the spacecraft rises above the $y$-$z$ plane. Because the sublimation stops at the beginning of the deflection operations, the efficiency of the deflection, in the case of the natural formations, is strongly affected by the position along the orbit at which the sublimation starts. This is consistent with the results presented in \cite{colombo2009b}.



\section{Effect of Eccentricity}
One may argue that the method is effective only on asteroids relatively close to the Sun as the solar collectors need to power the laser. Indeed if the asteroid has an aphelion far from the Sun the power can drop below the minimum required to sublimate the surface. Using the idea of the shaped orbits, one can try to apply the laser concept to asteroids with an increasing aphelion from 1 AU to 2 AU and with a decreasing perihelion from 1 AU to 0.5 AU. The assumption is that the Earth is moving on a circular planar orbit and the asteroid on a planar elliptic orbit. The impact parameter is computed at one of the two intersections with the orbit of the Earth and the deflection action starts at the perihelion of the orbit of the asteroid.

\begin{figure}[!htb]
\begin{center}
  \includegraphics[width=0.9\textwidth]{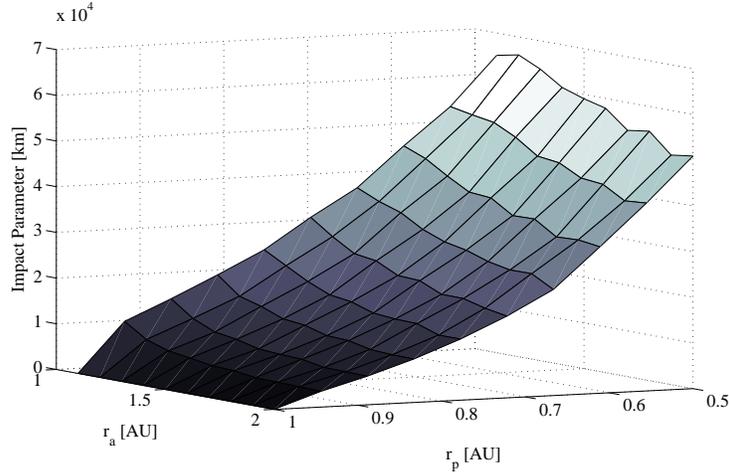}\\
  \caption{Impact parameter as a function of the radius of perihelion $r_P$ and aphelion $r_A$ of the orbit of the asteroid.}\label{fig:impact_rp_ra}
\end{center}
\end{figure}

Figure \ref{fig:impact_rp_ra} shows the achievable impact parameter as a function of radius of the aphelion and perihelion for 9 years of warning time, $C_r=5000$ and a 20 m diameter collector. For comparison with the case of Apophis one can notice that the achievable impact parameter is substantially high for highly elliptical asteroids. There are two good reasons for that. One is that the thrust is applied mainly at the pericentre of the orbit but for highly elliptical orbits a variation of velocity at the pericentre produces a much higher change of the semi-major axis than for low eccentric orbits. The other is that the orbit of the asteroid has a much steeper intersection with the Earth's orbit and therefore a small variation of the arrival time generates a greater impact parameter.

If one sticks to the hypothesis used above for the contamination, even in the case of natural orbits the spacecraft will experience no contamination as they fly above the plume when the sublimation is minimal or null. Therefore, the laser ablation seems to be effective even for high elliptical asteroids with high aphelion.

\section{Conclusion}
This paper presented the multidisciplinary design of a formation of spacecraft equipped with solar pumped laser for the deflection of asteroids.

The paper demonstrated that the use of multiple spacecraft is an optimal solution to maximise the deflection while minimizing the mass of the overall system. In fact as the diameter of the primary mirror increases the radiator and laser mass increases up to the point at which the mass of a single spacecraft exceeds the total mass of two or more spacecraft of smaller size. This is a very important point that is in favour of the use of a formation instead of a single large spacecraft. A formation, or fractionated system, has the further advantage of increasing redundancy and scalability as for a bigger asteroid the solution is simply to increase the number of spacecraft. The sizing of the spacecraft was based on a simple model in which the mass of the main bus is considered constant and the propellant mass is not optimised. These are two limiting assumptions that cause an overestimation of the mass for small systems. At the same time the deployment and thermal control systems are assumed to be scalable within the range of variability of the design parameters. Looking at present technology, this assumption can correspond to an underestimation of the mass for large systems. The efficiency of the laser and solar cells are at the upper limit of what is currently achievable in a lab environment. Although this is an optimistic assumption, current developments are progressing towards those limits independently of the deflection of asteroids. It is therefore reasonable to expect the system efficiencies presented in this paper in the near future. The paper also analyzed the control of the spacecraft in the vicinity of the asteroid and showed that with minimal control and propellant consumption the spacecraft can be maintained in their desired formation orbits.

Finally it was demonstrated that the laser ablation concept based on solar power is applicable also to high eccentric orbits (deep crossers) with even better performance with respect to the shallow crosser case. In fact, for deep crossers the deflection action is maximal where most effective, i.e., around the perihelion, and the steep intersection between orbit of the Earth and orbit of the asteroid amplifies the deflection effect.

\section{Acknowledgements}
This research was partially supported by the ESA/Ariadna Study Grant AO/1-5387/07/NL/CB \citep{ariadna}. The authors would like to thank Dr.\ Leopold Summerer of the ESA Advanced Concepts Team for his support.


\bibliographystyle{elsarticle-harv}
\bibliography{bee_design}   

\end{document}